\documentclass[10pt]{article}
\usepackage{flafter,amsmath,amssymb,latexsym,psfrag,graphicx,color}
\usepackage[margin=2.5cm]{geometry}
\usepackage[numbers,sort&compress]{natbib}

\newtheorem{theorem}{Theorem}[section]
\newtheorem{defn}[theorem]{Definition}
\newtheorem{lemma}[theorem]{Lemma}
\newtheorem{remark}[theorem]{Remark}

\newtheorem{cor}[theorem]{Corollary}
\newtheorem{prop}[theorem]{Proposition}

\allowdisplaybreaks

\makeatletter

\@addtoreset{figure}{section}
\makeatother

\begin{document}
\setlength\arraycolsep{2pt}
\date{\today}

\title{Estimation of the heat conducted by a cluster of small cavities and characterization of the equivalent heat conduction}
\author{Mourad Sini$^1$, Haibing Wang$^{2,\,3}$
\\$^1$RICAM, Austrian Academy of Sciences, A-4040, Linz, Austria \qquad\qquad\quad\quad
\\E-mail: mourad.sini@oeaw.ac.at
\\$^2$School of Mathematics, Southeast University, Nanjing 210096, P.R. China\;
\\$^3$Shing-Tung Yau Center, Southeast University, Nanjing 210096, P.R. China
\\E-mail: hbwang@seu.edu.cn
}

\maketitle
\begin{abstract}

We estimate the heat conducted by a cluster of many small cavities. We show that the dominating heat is a sum, over the number of the cavities, of the heats generated by each cavity after interacting with each other. This interaction is described through densities computable as solutions of a close, and invertible, system of time domain integral equations of a second kind. As an application of these expansions, we derive the effective heat conductivity which generates approximately the same heat as the cluster of cavities, distributed in a $3$D bounded domain, with explicit error estimates in terms of that cluster. At the analysis level, we use time domain integral equations. Doing that, we have two choices. First, we can favor the space variable by reducing the heat potentials to the ones related to the Laplace operator (avoiding Laplace transform). Second, we can favor the time variable by reducing the representation to the Abel integral operator. As the model under investigation has time-independent parameters, we follow here the first approach.



\end{abstract}

\section{Introduction}
\setcounter{equation}{0}

Let $D$ be a bounded domain in $\mathbb R^3$ and consider the following initial-boundary value problem:
\begin{equation}\label{bvp}
\begin{cases}
(\partial_t - \Delta)u=0 & \textrm{ in } (\mathbb R^3\setminus\overline D)_T, \\
u=f & \textrm{ on } (\partial D)_T, \\
u=0 & \textrm{ at } t=0.
\end{cases}
\end{equation}
For simplicity of notations, here and throughout this paper, we denote $X\times (0,\,T)$ and $\partial X\times (0,\,T)$ by $X_T$ and $(\partial X)_T$, respectively, where $X$ is a domain in $\mathbb R^3$ and $\partial X$ denotes its boundary. We also set $T_\varepsilon:=T/\varepsilon^2$ where $\varepsilon$ is a positive parameter.
The model \eqref{bvp} has many important applications in sciences and engineering \cite{N-W2013, N-W2015, N-W2017}. For example, in active thermography, $D$ is regarded as a cavity embedded in the background medium.

We assume in \eqref{bvp} that the compatibility condition $f(x,\,0)=u(x,\,0)$ holds on $\partial D$. To ensure the uniqueness of solutions to \eqref{bvp}, we require that the solution $u(x,\,t)$ satisfies the growth condition
\begin{equation}\label{growth}
|u(x,\,t)|\leq C_0 \exp(b|x|^2) \qquad \textrm{as}\; |x| \to +\infty
\end{equation}
for some positive constants $C_0$ and $b < (4T)^{-1}$. Under certain regularity assumptions on $D$ and $f$, the unique solvability of the problem \eqref{bvp} can be proved; see for instance \cite{Pao1992}.

In this work, we consider the case where $D$ is given by a union of small cavities, i.e. $D:=\cup^M_{j=1} D_j$, with the maximum radius of $D_j$ of small order $\varepsilon$. Our goal is to estimate the heat generated by such a cluster in terms of the size and the number of the cavities. To give sense to the boundary conditions on the surfaces of the small cavities, we take a source defined, at least, in a domain containing the whole cluster and then the boundary conditions on each cavity's surface is the trace of this source on it.

We could consider more general sources in the model $(\ref{bvp})$, namely
\begin{equation}\label{bvp-2}
\begin{cases}
(\partial_t - \Delta)w=F & \textrm{ in } (\mathbb R^3\setminus\overline D)_T, \\
w=G & \textrm{ on } (\partial D)_T, \\
w=H & \textrm{ at } t=0
\end{cases}
\end{equation}
with the needed compatibility conditions for the sources $G$ and $H$, i.e. $G(\cdot,\,0)=H(\cdot)$ on $\partial D$. However, setting $v$ to be the solution in the absence of the cavities, by extending first $F$ and $G$ to the whole space $\mathbb{R}^3$, we see that $u:=w-v$ solves $(\ref{bvp})$ with the boundary source $f:=G-v|_{(\partial D)_T}$.
A particular source of heat that is used in practice is given by
\begin{equation}\label{source}
f(x,\,t):=f(x,\,t;\,z^*)=\Phi(x,\,t;\,z^*,\,0),
\end{equation} where
\begin{equation}\label{fund}
\Phi(x,\, t;\, y,\, \tau):=
\begin{cases}
\displaystyle\frac{1}{[4\pi (t-\tau)]^{3/2}} \exp\left( -\frac{| x-y |^2}{4(t-\tau)} \right), & t>\tau,\\
0, & t\leq \tau
\end{cases}
\end{equation}
is the fundamental solution of the heat operator $\partial_t - \Delta$ for the three-dimensional spatial space and $z^*$ is the source point which is located away from the cluster $\cup^M_{j=1} D_j$. This source is initially a Dirac source supported on the source point $z^*$, i.e. $f(x,\,0)=\delta (z^*)$. With such sources, we see from (\ref{bvp}), the $u_f:=f-u$ satisfies the problem
\begin{equation}\label{bvp-f}
\begin{cases}
(\partial_t - \Delta)u_f=0 & \textrm{ in } (\mathbb R^3\setminus\overline D)_T, \\
u_f=0 & \textrm{ on } (\partial D)_T, \\
u_f=\delta(z^*) & \textrm{ at } t=0
\end{cases}
\end{equation}
with the growth condition as in (\ref{growth}).
\bigskip

Let $B_1,\,B_2,\,\cdots,\,B_M$ be $M$ open, bounded and simply connected domains in $\mathbb R^3$ with $C^2$-boundaries containing the origin. Assume that the Lipschitz constants of $B_j,\,j=1,\,2,\,\cdots,\,M$ are uniformly bounded. Set $D_j := \varepsilon B_j + z_j$ to be small cavities characterized by the parameter $\varepsilon >0$ and the locations $z_j\in\mathbb R^3,\,j=1,\,2,\,\cdots,\,M$.

\begin{defn}\label{defn}
Let $D$ be a union of the small cavities, i.e., $D:=\cup_{j=1}^{j=M}D_j$. We set
\begin{enumerate}
\item $ a$ as the maximum among the diameters of the small cavities, i.e.,
\begin{equation}\label{defn_a}
a:=\max_{1\leq j \leq M}\,\mathrm{diam}(D_j)=\varepsilon \, \max_{1\leq j \leq M}\,\mathrm{diam}(B_j);
\end{equation}

\item $d$ as the minimum distance between the small cavities, i.e.,
\begin{equation}\label{defn_d}
d:=\min_{\substack{1\leq i,j \leq M\\i\neq j}}\,d_{ij}, \quad d_{ij}:= \mathrm{dist}(D_i,\,D_j).
\end{equation}
\end{enumerate}
\end{defn}

We are interested in regimes where
\begin{equation}\label{regimes}
M\sim a^{-s}\; \mbox{ and }\; d\sim a^{\beta}
\end{equation}
with nonnegative real numbers $s$ and $\beta$.
\bigskip

Our first result is the following approximation property:

\begin{theorem}\label{Main-asymptotic}
Let $D_i, i=1, ..., M,$ and the source $f$ be as described above. Under the condition
\begin{equation}\label{cluster-condition}
a\; \max_{1\leq i \leq M}\sum_{j\neq i}d^{-2}_{ij}<1,
\end{equation}
which means that $1-2\beta-\frac{s}{3}\geq 0$, the heat conducted by this cluster is approximately estimated, for $x$ away from $\cup_{j=1}^{j=M}D_j$ and $t>0$, as
\begin{equation}
u(x,\,t) = \sum_{i=1}^M C_i \int_0^t \Phi(x,\,t;\,z_i,\,\tau)\, \alpha_i(\tau)\, d\tau +
O(a^{3-3\beta-s} |\ln a|) +O(a^{2-s}) \quad \textrm{as}\; a \to 0,
\end{equation}
where each constant $C_i$ is the capacitance of the cavity $D_i,\,i=1,\,\cdots,\, M$, defined as  $
 C_i:=\int_{\partial D_i}\sigma_i(x)\,ds(x)$
with $\sigma_i$ as the unique solution of the integral equation $\int_{\partial D_i}\frac{\sigma_i(y)}{4\pi \vert x-y \vert}\,ds(y) =1, \; x\in \partial D_i,$
and $\left\{ \alpha_i(t) \right\}_{i=1}^M$ is the unique solution of the linear system
\begin{equation}\label{F-L-system}
\alpha_i(t) + \sum_{\substack{j=1\\j\neq i}}^M \int_0^t C_j \,\Phi(z_i,\,t;\,z_j,\,\tau)\,\alpha_j(\tau)\, d\tau = \Phi(z_i,\,t;\,z^*,\,0),\quad t\in (0,\,T),\quad i=1,\,2,\,\cdots,\,M,
\end{equation}
which is invertible from $L^2(0,\,T)$ to itself.
\end{theorem}

We can show that $\max\limits_{1\leq i \leq M}\sum\limits_{j\neq i}d^{-2}_{ij}\leq d^{-2}M^{\frac{1}{3}}$, see \cite{BS2018} for instance, then (\ref{cluster-condition}) is satisfied if $a\, d^{-2} M^{\frac{1}{3}}<1$, and from (\ref{regimes}), it implies that $1-2\beta -\frac{s}{3}\geq 0$.

\bigskip
As an application of such results, we derive the effective heat conductivity distribution that can produce the same heat as the cluster above.
To show this, let $\Omega$ be a bounded domain containing the cavities $D_j,\,j=1,\,2,\,\cdots,\,M$. We divide $\Omega$ into $[a^{-1}]$ subdomains
$\Omega_j,\,j=1,\,2,\,\cdots,\,[a^{-1}]$, periodically arranged for instance\footnote{The periodicity is actually not needed. We assume it only for simplicity of exposition.},
such that the $\Omega_j$'s are disjoint and of a volume $a$. Here, we denote by $[x]$ the unique integer $n$ such that $n\leq x < n+1$, i.e. $n$ is the floor number.
Each subdomain $\Omega_j$ contains one single hole. Such a distribution obeys the condition (\ref{cluster-condition}).
Indeed, from the estimate
$\max_{1\leq i \leq M}\sum_{j\neq i}d^{-2}_{ij}\leq d^{-2} M^{\frac{1}{3}}$
and as $M=[a^{-1}]\leq a^{-1}$, i.e. $s=1$, we have
$a\, \max_{1\leq i \leq M}\sum_{j\neq i}d^{-2}_{ij}\leq a^{\frac{2}{3}} d^{-2}=\left(a^{\frac{1}{3}}d^{-1}\right)^2.$
Taking $d=d_0a^{\frac{1}{3}}$, i.e. $\beta=\frac{1}{3}$, we see $a\, \max_{1\leq i \leq M}\sum_{j\neq i}d^{-2}_{ij}\leq \frac{1}{d^2_0}<1$ if $d_0>1.$ This last condition is obviously satisfied according to the distribution described above.

Now, we state our second main result.
\begin{theorem}\label{Main-effective}
Let $\Omega$ be a bounded and Lipschitz domain in $\mathbb{R}^3$. We distribute the cavities as described above. Then for any $t>0$ and $x$ away from $\Omega$, we have the approximation:
\begin{equation}
u(x,\,t)-W(x,\,t)=O(a^{\frac{1}{3}})\; \mbox { as } a\ll 1,
\end{equation}
where $W$ is the unique solution of the problem
\begin{equation}\label{effective-medium-W}
\begin{cases}
(\partial_t - \Delta + \overline C\chi_\Omega)W = \overline C \chi_\Omega \Phi(x,\,t;\,z^*,\,0) & \mathrm{in}\;\mathbb R^3\times (0,\,T), \\
W(x,\,0)=0 & \mathrm{in}\; \mathbb R^3, \\
|W(x,\,t)|\leq C_0 \exp(b|x|^2) & \mathrm{as}\; |x|\to +\infty
\end{cases}
\end{equation}
for some positive constants $C_0$ and $b < (4T)^{-1}$. Here $\overline C$ is the capacitance of the unscaled domains $B_m$'s (that are assumed to be the same). Finally $\chi_\Omega$ is the characteristic function of the domain $\Omega$.
\end{theorem}

As a corollary, we deduce the following result. Let $\sigma$ be the unique solution of the problem
$-\Delta \sigma +\overline C \sigma =0 $ in $\Omega$ and $\sigma =1$ on $\partial \Omega$.
As $\overline C$ is positive in $\Omega$, by the maximum principle, the unique solution of this problem is positive in $\Omega$. We extend $\sigma$ from $\Omega$ to $\mathbb{R}^3$ simply by $1$ in $\mathbb{R}^3\setminus{\overline \Omega}$. By a simple change of variable $\tilde{W}:=\sigma^{-1}\; W$, we see that $\tilde{W}$ satisfies the problem
\begin{equation}\label{effective-medium-tilde-W}
\begin{cases}
\sigma\,\partial_t \tilde W - \sigma^{-1}\nabla \cdot \sigma^2 \nabla \tilde{W} = \overline C \chi_\Omega \Phi(x,\,t;\,z^*,\,0) & \mathrm{in}\;\mathbb R^3\times (0,\,T), \\
\tilde{W}(x,\,0)=0 & \mathrm{in}\; \mathbb R^3, \\
|\tilde{W}(x,\,t)|\leq C_0 \exp(b|x|^2) & \mathrm{as}\; |x|\to +\infty.
\end{cases}
\end{equation}
Observe that the first equation in (\ref{effective-medium-tilde-W}) can be written as $\rho \, c\, \partial_t \tilde W - \nabla \cdot \gamma \nabla \tilde{W} = F $ where the density $\rho$ and the heat capacity $c$ are such that $\rho \, c=\sigma^2$ and the heat conductivity $\gamma$ is give by $\gamma:=\sigma^2$.

As $\tilde{W}=W$ in $\mathbb{R}^3\setminus{\overline \Omega}$, from Theorem \ref{Main-effective} we deduce the following result:

\begin{cor} \label{Main-effective-cor}
Let $\Omega$ as described in Theorem \ref{Main-effective}. Then for any $t>0$ and $x$ away from $\Omega$, we have the approximation:
\begin{equation}
u(x,\,t)-\tilde{W}(x,\,t)=O(a^{\frac{1}{3}})\; \mbox { as } a \ll 1,
\end{equation}
where $\tilde{W}$ is the unique solution of the problem (\ref{effective-medium-tilde-W}) with $\sigma$ being the unique solution of the problem
\begin{equation}\label{sigma}
\begin{cases}
-\Delta \sigma +\overline C \sigma =0 & \mathrm{ in }\; \Omega,\\
 \sigma =1 & \mathrm{ on }\; \partial \Omega.
\end{cases}
\end{equation}
\end{cor}

The results provided in Theorem \ref{Main-asymptotic}, Theorem \ref{Main-effective} and Corollary \ref{Main-effective-cor} are given for boundary sources $f(x,\,t)$ as point sources initially located at $z^*$.
Actually, these results are valid for any source $f$ in the Sobolev space $H_0^1((0,\,T);\, W^{1,\,\infty}(\mathbb{R}^3))$.

\medskip
Next, observe that $G_C(x,\,t;\, z^*,\, 0):= \Phi(x,\,t;\,z^*,\,0)-W(x,\,t)$ satisfies
\begin{equation}\label{effective-medium-G}
\begin{cases}
(\partial_t - \Delta + \overline C\chi_\Omega)G_C = \delta_x(z^*)\delta_t(0) & \mathrm{in}\;\mathbb R^3\times (0,\,T), \\
G_C(x,\,0)=0 & \mathrm{in}\; \mathbb R^3, \\
|G_C(x,\,t)|\leq C_0 \exp(b|x|^2) & \mathrm{as}\; |x|\to +\infty.
\end{cases}
\end{equation}
Now, using the form $f(x,\,t):=\Phi(x,\,t;\,z^*,\,0)$ in (\ref{bvp}), we see that $G_D(x,\,t;\,z^*,\,0):= \Phi(x,\,t;\,z^*,\,0)-u(x,\,t)$ satisfies
\begin{equation}\label{bvp-GD}
\begin{cases}
(\partial_t - \Delta)G_D=\delta_x(z^*)\delta_t(0) & \textrm{ in } (\mathbb R^3\setminus\overline D)_T, \\
G_D=0 & \textrm{ on } (\partial D)_T, \\
G_D=0 & \textrm{ at } t=0.
\end{cases}
\end{equation}
From Theorem \ref{Main-effective}, we see that, for $t>0$ and $ x \in \mathbb{R}^3\setminus{\overline{\Omega}}$, we have
$$
G_D(x,\,t;\,z^*,\,0)-G_C(x,\,t;\,z^*,\,0)=O(a^{\frac{1}{3}})\; \mbox{ as } a \ll 1.
$$
As $z^*$ is arbitrary in $\mathbb{R}^3\setminus{\overline \Omega}$, this approximation implies that for any given initial source function $H: =H(x)$ compactly supported in $\mathbb{R}^3\setminus{\overline \Omega}$, then for $x\in \mathbb{R}^3\setminus{\overline \Omega}$ and $t>0$, we have
$$
(G_D \ast H) (x,\,t)-(G_C \ast H)(x,\,t)=O(a^{\frac{1}{3}})\; \mbox{ as } a \ll 1
$$
where $G_D \ast H$ and $G_C\ast H$ stand for the convolutions in the space variable of the Green kernels $G_D$ and $G_C$, respectively with $H$. Hence, $G_D \ast H$ is the solution of (\ref{bvp-2}) with $H$ as the initial data, with $F=0$ and $G=0$, i.e.
\begin{equation}\label{bvp-3}
\begin{cases}
(\partial_t - \Delta)(G_D \ast H)=0 & \textrm{ in } (\mathbb R^3\setminus\overline D)_T, \\
(G_D \ast H)=0 & \textrm{ on } (\partial D)_T, \\
(G_D \ast H)=H & \textrm{ at } t=0
\end{cases}
\end{equation}
and $G_C\ast H$ is the one of the problem
\begin{equation}\label{effective-medium-G-f}
\begin{cases}
(\partial_t - \Delta + \overline C\chi_\Omega)(G_C\ast H) = 0 & \mathrm{in}\;\mathbb R^3\times (0,\,T), \\
(G_C\ast H)(x,\,0)=H & \mathrm{in}\; \mathbb R^3, \\
|(G_C\ast H)(x,\,t)|\leq C_0 \exp(b|x|^2) & \mathrm{as}\; |x|\to +\infty.
\end{cases}
\end{equation}

Few remarks are in order:

\begin{enumerate}
\item The observation described above is not a surprise. That means, estimating the heat, generated by the cluster $\cup^M_{i=1} D_i$, starting from an initial heat supported on the source points $z^*$ arbitrary located outside $\Omega$ is enough to estimate the heat starting from any source supported in $\mathbb{R}^3\setminus{\overline{\Omega}}$.

\item The condition that the initial source $H$ is supported outside $\Omega$ is not a restriction. Indeed, see (\ref{bvp-3}), as $H=0$ on $\cup^M_{i=1} \partial D_i$ and the set of cavities is densely distributed in $\Omega$, then necessarily we should assume that $H=0$ in $\Omega$.

\item The condition (\ref{cluster-condition}) can be relaxed as follows. Recall that this condition translates as $1-2\beta -\frac{s}{3}\geq 0$. This condition is enough to derive the effective conductivity where we need $s=1$ and $\beta =\frac{1}{3}$. However, performing the analysis done in this work more carefully, we can handle larger regimes where $s>1$ and $\beta >\frac{1}{3}$.

\item In Theorem \ref{Main-effective}, we divided $\Omega$ into a family $\Omega_m$'s which contain a single cavity each. Actually, we can put arbitrary number of cavities in each $\Omega_m$. This would be translated by the appearance of the local distribution density, we denote by $K:=K(z)$, in (\ref{effective-medium-W}), replacing $\overline{C}$ by $\overline{C} K$.

\item In Corollary \ref{Main-effective-cor}, we have seen that the cluster behaves as a conductive heat medium modeled by a heat conductivity $\gamma:=\sigma^2$.
This is possible because such cavities, modeled by Dirichlet boundary conditions, are actually resonating and hence enhance the generated heats.
This is translated by the fact that the additive potential $\overline{C} \chi_\Omega$ appears with a positive sign and hence the generated conductivity $\sigma$
via the boundary value problem (\ref{sigma}) is positive. This phenomenon can also occur for other resonating particles as the electromagnetic nanoparticles enjoying
balanced contrast/size ratios. However, this is not always the case for other boundary conditions, as the impedance boundary conditions, for which the sign of the corresponding coefficient $\overline C$ can be negative.
In this case, the corresponding boundary value problem (\ref{sigma}) does not enjoy positivity for its solutions and hence the heat conductivity $\sigma$ might not be generated.

\item In Theorem \ref{Main-effective}, we distributed the cavities in $3$D domains. It is natural, and interesting in applications, to consider distributions in $2$D surfaces or $1$D curves and their different superpositions. The corresponding results for these situations will be reported elsewhere.

\end{enumerate}

The estimation of the fields generated by a cluster of small particles (of different kinds) is well developed in the literature for the elliptic models, both stationary or nonstationary case, see \cite{B-L-P:1978, DC-FM:Book:BB1997, J-K-O:1994, M-K:2006} for periodic distributions and \cite{A-C-K-S-15, A-K:2007, Bendali-et.al:2015, BS2018, C-S2014, Martin:2006, M-M-N:Springerbook:2013, RAMM:2007} for nonperiodic distributions. This list is of course by no means exhaustive. The situation is much less clear for time domain models, as those related to parabolic, Schr\"odinger or hyperbolic equations, unless for periodic media, see for instance \cite{DC-FM:Book:BB1997, M-K:2006}. Nevertheless, for the particular situations of point-like particles, there are several works published in the framework of singular perturbations, see for instance \cite{ABD:95, DFT:96, DFT:98}.

Motivated by applications in mathematical imaging, asymptotic expansions of the heat generated by single, or well separated, small particles are derived in \cite{AIKK2005, ARR:2018} based on energy methods or Laplace transform. Here, we avoid using Laplace transform and rather, we use time domain integral equations. One of the basic arguments we used here can be explained as follows. Based on integral representations, via the single-layer (or double-layer) heat potentials, we have two choices. Either, we favor the space variable, i.e. reduce the single-layer heat potential to the one related to the Laplace operator, or the time variable, i.e. reduce the representation to the Abel integral.
In this paper, since the cavities are fixed and not moving in time, we follow the first approach as it allows us to extract in a straightforward way the dominant part of the heat generated by the cluster. However, the other approach is also interesting as it might allow moving cavities (for instance allowing their centers to be time dependent).

The rest of the paper is organized as follows. In Section \ref{preliminaries}, as preliminaries, we recall some properties on the integral equation method for the heat equation, mainly the single-layer heat potential. In Section \ref{single-cavity}, we provide the analysis for the case of a single cavity to describe the main steps of our approach. In Section \ref{multiple-cavities}, we provide the proof of Theorem \ref{Main-asymptotic} and in Section \ref{effective-medium}, the one of Theorem \ref{Main-effective}.

\section{Preliminaries}\label{preliminaries}
\setcounter{equation}{0}

Throughout the paper, we use the notation \lq\lq$\lesssim$\rq\rq{} to denote \lq\lq$\leq$\rq\rq{} with its right-hand side multiplied by a generic positive constant.
We recall some known properties of the single-layer heat operator. We have the following classical singularity estimates for the fundamental solution $\Phi(x,\,t;\,y,\,\tau)$:
\begin{eqnarray}
&&|\Phi(x,\,t;\,y,\,\tau)| \lesssim  \frac{1}{(t-\tau)^\mu} \,\frac{1}{|x-y|^{3-2\mu}},\quad \mu<\frac{3}{2}, \label{fun_est1}\\
&&|\partial_{x_j}\Phi(x,\,t;\,y,\,\tau)| \lesssim \frac{1}{(t-\tau)^\gamma} \,\frac{1}{|x-y|^{4-2\gamma}},\quad \gamma<\frac{5}{2},\;j=1,\,2,\,3, \label{fun_est2}\\
&&|D_t^{1/2}\Phi(x,\,t;\,y,\,\tau)| \lesssim \frac{1}{(t-\tau)^\gamma} \,\frac{1}{|x-y|^{4-2\gamma}},\quad \gamma<\frac{5}{2} \label{fun_est3}
\end{eqnarray}
for $0\leq \tau <t\leq T$ and $x,\,y\in\mathbb R^3$ with $x\neq y$. Here $D_t^{1/2}$ denotes the time derivative of order $1/2$ defined by
\begin{equation*}
D_t^{1/2}h(t):= \frac{1}{\sqrt \pi}\int_0^t \frac{h^\prime(s)}{\sqrt{t-s}}\,ds.
\end{equation*}

Define the single-layer heat potential $\mathcal S_{\Gamma_T}$ by
\begin{equation}\label{defn_S}
\mathcal S_{\Gamma_T}[\varphi](x,\,t):=\int_0^t\int_\Gamma \Phi(x,\,t;\,y,\,\tau)\,\varphi(y,\,\tau) \,ds(y)d\tau,\quad (x,\,t)\in \Gamma\times(0,\,T),\\
\end{equation}
where $\Gamma$ is the boundary of a bounded domain in $\mathbb R^3$. To look for the weak solution of \eqref{bvp} by the single-layer potential, we introduce the anisotropic Sobolev space $H^{1,1/2}(\Gamma_T)$.
It is a closure of the space
\begin{equation*}
\left\{w: \,w=v|_{\Gamma_T},\,v\in C_0^\infty\left(\mathbb R^3 \times (0,\,+\infty)\right)\right\}
\end{equation*}
with respect to the norm
\begin{equation}\label{defn_H1}
\| w\|^2_{H^{1,1/2}(\Gamma_T)}:=\int_0^T \int_\Gamma \left[|\nabla_{\mathrm{tan}} w(x,\,t)|^2 + w^2(x,\,t) + \left(D_t^{1/2}w(x,\,t)\right)^2\right]\, ds(x)dt ,
\end{equation}
where $\nabla_{\mathrm{tan}} w$ is the tangential gradient of $w$ on $\Gamma$ defined by $
\nabla_{\mathrm{tan}} w := \nabla w - (\partial_\nu w)\nu $
with $\nu$ being the outward unit normal vector to $\Gamma$. Then the invertibility of $\mathcal S_{\Gamma_T}$ can be stated as follows \cite{Brown1987, Brown1989}.

\begin{lemma}\label{S_inv1} Assume $\Gamma$ to be Lipschitz regular.
The operator $\mathcal S_{\Gamma_T}:\, L^2(\Gamma_T)\to H^{1,1/2}(\Gamma_T)$ is invertible with a bounded inverse.
\end{lemma}

To investigate the classical solution of \eqref{bvp} and consider its regularity, we define the anisotropic H\"older spaces of functions
\begin{equation*}
C_\lambda^{k,\alpha}(\Gamma_T)=\left\{ u \big|\, \|u\|_\lambda^{(k,\,\alpha)}(\Gamma_T) < +\infty \right\},\quad k=0,\,1,\;\alpha>0,\;\lambda>0,
\end{equation*}
where
\begin{eqnarray*}
\|u\|_\lambda^{(0,\,\alpha)}(\Gamma_T)&:=& \sup_{(x,\,t)\in\Gamma_T} \frac{|u(x,\,t)|}{\exp\lambda t} + \sup_{\substack{(x,\,t),\,(x+\Delta x,\,t+\Delta t)\in\Gamma_T\\|(\Delta x,\,\Delta t)|\neq 0}} \frac{|u(x+\Delta x,\,t) - u(x,\,t)|}{|\Delta x|^\alpha [\exp\lambda t + \exp\lambda(t+\Delta t)]}, \\
\|u\|_\lambda^{(1,\,\alpha)}(\Gamma_T)&:=& \sup_{(x,\,t)\in\Gamma_T} \frac{|u(x,\,t)|}{\exp\lambda t} + \sup_{\substack{(x,\,t),\,(x,\,t+\Delta t)\in\Gamma_T\\|\Delta t|\neq 0}} \frac{|u(x,\,t+\Delta t) - u(x,\,t)|}{|\Delta t|^{(1+\alpha)/2} [\exp\lambda t + \exp\lambda(t+\Delta t)]}\\
 && \qquad\qquad\qquad\quad\; + \sum_{i=1}^3 \|\partial_{x_i} u \|_\lambda^{(0,\,\alpha)}(\Gamma_T).
\end{eqnarray*}
Define $\mathop{C}\limits_{\circ}{}_\lambda^{k,\alpha}(\Gamma_T)$ by
\begin{equation*}
\mathop{C}\limits_{\circ}{}_\lambda^{k,\alpha}(\Gamma_T):= \left\{u\in C_\lambda^{k,\alpha}(\Gamma_T)\,\big|\, u(x,\,0)=0 \right\}.
\end{equation*}
Then we have the following result \cite{Baderko1992, Baderko1997}:
\begin{lemma}\label{S_inv2} Assume $\Gamma$ to be of class $C^2$.
The operator $\mathcal S_{\Gamma_T}:\, \mathop{C}\limits_{\circ}{}_\lambda^{0,\alpha}(\Gamma_T)\to \mathop{C}\limits_{\circ}{}_\lambda^{1,\alpha}(\Gamma_T)$ is invertible with a bounded inverse.
\end{lemma}
Observe that by Sobolev embedding, if $f \in H_0^1((0,\,T);\, W^{1,\,\infty}(\mathbb{R}^3))$ then $f|_{\Gamma_T} \in \mathop{C}\limits_{\circ}{}_\lambda^{0,\frac{1}{2}}(\Gamma_T)$.
Based on Lemma \ref{S_inv1}, we see that the problem (\ref{bvp}) is well-posed for sources $f$ belonging to $L^2(\Gamma_T)$. In addition, based on Lemma \ref{S_inv2},
the problem enjoys regularity properties. A particular property that we need is that if $f(\cdot,\,0)=0$, then also the density $\phi$ satisfies the same property, i.e. $\phi(\cdot,\,0)=0$.
Actually, in our analysis we need this property and this is the only place where we need the $C^2$ regularity of $\Gamma$. Otherwise, only the Lipschitz smoothness of $\Gamma$ is needed.

\medskip
Let $D=\varepsilon B + z$, where $B$ is a bounded and simply connected domain in $\mathbb R^3$ with Lipschitz boundary containing the origin. For any functions $\varphi$ and $\psi$ defined on $(\partial D)_T$ and $(\partial B)_{T_\varepsilon}$, respectively, we use the notations
\begin{equation*}
\hat\varphi(\eta,\,\tilde \tau)=\varphi^\wedge(\eta,\,\tilde \tau):=\varphi(\varepsilon\eta+z,\,\varepsilon^2\tilde \tau),\quad \check\psi(x,\,t)=\psi^\vee(x,\,t):=\psi\left(\frac{x-z}{\varepsilon},\,\frac{t}{\varepsilon^2}\right)
\end{equation*}
for $(x,\,t)\in (\partial D)_T$ and $(\eta,\,\tilde\tau)\in (\partial B)_{T_\varepsilon}$. Then we have the following lemmas.

\begin{lemma}\label{lemma1}
Suppose $0<\varepsilon\leq 1$ and $D=\varepsilon B + z \subset \mathbb R^3$. Then, for $\varphi\in H^{1,1/2}((\partial D)_T)$ and $\psi\in L^2((\partial D)_T)$, we have
\begin{equation}\label{scalepsi}
\|\psi\|_{L^2((\partial D)_T)} = \varepsilon^2 \|\hat\psi\|_{L^2((\partial B)_{T_\varepsilon})}
\end{equation}
and
\begin{equation}\label{scalephi}
\varepsilon^2 \| \hat \varphi\|_{H^{1,1/2}((\partial B)_{T_\varepsilon})} \leq \| \varphi\|_{H^{1,1/2}((\partial D)_T)} \leq \varepsilon \| \hat \varphi\|_{H^{1,1/2}((\partial B)_{T_\varepsilon})}.
\end{equation}
\end{lemma}

{\bf Proof.} Let $x=\varepsilon \xi + z$ and $t=\varepsilon^2 \tilde t$. Then for $\psi\in L^2((\partial D)_T)$ we have
\begin{equation*}
\|\psi\|^2_{L^2((\partial D)_T)}= \int_0^T \int_{\partial D} |\psi(x,\,t)|^2\,ds(x)dt
= \varepsilon^4 \int_0^{T_\varepsilon}\int_{\partial B} |\psi(\varepsilon \xi + z,\,\varepsilon^2 \tilde t)|^2\,ds(\xi)d\tilde t
= \varepsilon^4 \|\hat\psi\|_{L^2((\partial B)_{T_\varepsilon})},
\end{equation*}
which leads to \eqref{scalepsi}.

For $\varphi\in H^{1,1/2}((\partial D)_T)$, we have
\begin{eqnarray*}
\| \varphi\|^2_{H^{1,1/2}((\partial D)_T)}&=& \int_0^T \int_{\partial D} \left[|\nabla_{\mathrm{tan}} \varphi(x,\,t)|^2 + \varphi^2(x,\,t) + \left(D_t^{1/2}\varphi(x,\,t)\right)^2\right] ds(x)dt\\
&=& \varepsilon^4 \int_0^{T_\varepsilon} \int_{\partial B} \left[\varepsilon^{-2}|\nabla_{\mathrm{tan}} \varphi(\varepsilon \xi + z,\,\varepsilon^2 \tilde t)|^2 + \varphi^2(\varepsilon \xi + z,\,\varepsilon^2 \tilde t)\right]\, ds(\xi)d\tilde t \\
& & + \varepsilon^2\int_0^{T_\varepsilon}\int_{\partial B} \left(D_{\tilde t}^{1/2}\varphi(\varepsilon \xi + z,\,\varepsilon^2\tilde t) \right)^2\, ds(\xi)d\tilde t,
\end{eqnarray*}
which implies \eqref{scalephi} by noticing that $\varepsilon^4 \leq \varepsilon^2$. \hfill $\Box$

\begin{lemma}\label{lemma2}
For $\varphi\in H^{1,1/2}((\partial D)_T)$ and $\psi\in L^2((\partial D)_T)$, we have
\begin{equation}\label{SD_1}
\mathcal S_{(\partial D)_T}[\psi] = \varepsilon (\mathcal S_{(\partial B)_{T_\varepsilon}}[\hat \psi])^\vee,
\end{equation}
\begin{equation}\label{SD_2}
\mathcal S_{(\partial D)_T}^{-1}[\varphi] = \varepsilon^{-1} (\mathcal S^{-1}_{(\partial B)_{T_\varepsilon}}[\hat \varphi])^\vee,
\end{equation}
and
\begin{equation}\label{SD_3}
\left \| \mathcal S_{(\partial D)_T}^{-1} \right\|_{\mathcal L \left(H^{1,1/2}((\partial D)_T),\,L^2((\partial D)_T)\right)} \leq \varepsilon^{-1}
\left \| \mathcal S_{(\partial B)_{T_\varepsilon}}^{-1} \right\|_{\mathcal L \left(H^{1,1/2}((\partial B)_{T_\varepsilon}),\,L^2((\partial B)_{T_\varepsilon})\right)}.
\end{equation}
\end{lemma}

{\bf Proof.} Let $x=\varepsilon \xi + z,\,y=\varepsilon \eta + z,\,t=\varepsilon^2 \tilde t,\,\tau=\varepsilon^2 \tilde \tau$. By direct calculations, we obtain
\begin{eqnarray*}
\mathcal S_{(\partial D)_T}[\psi](x,\,t) &=& \int_0^t \int_{\partial D} \displaystyle\frac{1}{[4\pi(t-\tau)]^{3/2}}\exp\left(-\frac{|x-y|^2}{4(t-\tau)}\right) \psi(y,\,\tau)\,ds(y)d\tau\\
&=& \int_0^{\tilde t} \int_{\partial B} \displaystyle\frac{1}{[4\pi\varepsilon^2(\tilde t-\tilde\tau)]^{3/2}}\exp\left(-\frac{|\xi-\eta|^2}{4(\tilde t-\tilde\tau)}\right) \psi(\varepsilon\eta+z,\,\varepsilon^2\tilde\tau)\, \varepsilon^4\,ds(\eta)d\tilde\tau\\
&=& \varepsilon \int_0^{\tilde t} \int_{\partial B} \displaystyle\frac{1}{[4\pi(\tilde t-\tilde\tau)]^{3/2}}\exp\left(-\frac{|\xi-\eta|^2}{4(\tilde t-\tilde\tau)}\right) \hat\psi(\eta,\,\tilde\tau)\,ds(\eta)d\tilde\tau\\
&=& \varepsilon \mathcal S_{(\partial B)_{T_\varepsilon}}[\hat\psi](\xi,\,\tilde t),
\end{eqnarray*}
which gives \eqref{SD_1}. Further, the identity \eqref{SD_2} follows from the following derivation:
\begin{equation*}
\mathcal S_{(\partial D)_T}[(\mathcal S^{-1}_{(\partial B)_{T_\varepsilon}}[\hat\phi])^\vee] = \varepsilon \left(\mathcal S_{(\partial B)_{T_\varepsilon}}\left[\mathcal S^{-1}_{(\partial B)_{T_\varepsilon}}[\hat\phi]\right]\right)^\vee = \varepsilon (\hat \varphi)^\vee = \varepsilon \varphi.
\end{equation*}

To show the estimate \eqref{SD_3}, we derive that
\begin{eqnarray*}
\left \| \mathcal S_{(\partial D)_T}^{-1} \right\|_{\mathcal L \left(H^{1,1/2}((\partial D)_T),\,L^2((\partial D)_T)\right)}
&:=& \sup\limits_{0\not=\varphi\in H^{1,1/2}((\partial D)_T)} \displaystyle \frac{ \| \mathcal S_{(\partial D)_T}^{-1}[\varphi] \|_{L^2((\partial D)_T)}}{\|\varphi\|_{H^{1,1/2}((\partial D)_T)}}\\
&\leq & \sup\limits_{0\not=\varphi\in H^{1,1/2}((\partial D)_T)} \displaystyle \frac{ \varepsilon^2\| (\mathcal S_{(\partial D)_T}^{-1}[\varphi])^\wedge \|_{L^2((\partial B)_{T_\varepsilon})}}{\varepsilon^2\|\hat\varphi\|_{H^{1,1/2}((\partial B)_{T_\varepsilon})}}\\
&= & \sup\limits_{0\not=\hat\varphi\in H^{1,1/2}((\partial B)_{T_\varepsilon})} \displaystyle \frac{ \| (\varepsilon^{-1}(\mathcal S_{(\partial B)_{T_\varepsilon}}^{-1}[\hat\varphi])^\vee)^\wedge \|_{L^2((\partial B)_{T_\varepsilon})}}{\|\hat\varphi\|_{H^{1,1/2}((\partial B)_{T_\varepsilon})}}\\
&= & \sup\limits_{0\not=\hat\varphi\in H^{1,1/2}((\partial B)_{T_\varepsilon})} \displaystyle \frac{\varepsilon^{-1} \| \mathcal S_{(\partial B)_{T_\varepsilon}}^{-1}[\hat\varphi] \|_{L^2((\partial B)_{T_\varepsilon})}}{\|\hat\varphi\|_{H^{1,1/2}((\partial B)_{T_\varepsilon})}}\\
&=&  \varepsilon^{-1}
\left \| \mathcal S_{(\partial B)_{T_\varepsilon}}^{-1} \right\|_{\mathcal L \left(H^{1,1/2}((\partial B)_{T_\varepsilon}),\,L^2((\partial B)_{T_\varepsilon})\right)}.
\end{eqnarray*}
Thus, the proof is complete. \hfill $\Box$

Finally, we need the following result \cite{A-N1989, Brown1987}:

\begin{lemma}\label{S_ep}
The single-layer potential operator $\mathcal S_{(\partial B)_{T_\varepsilon}}:\,L^2((\partial B)_{T_\varepsilon}) \to H^{1,1/2}((\partial B)_{T_\varepsilon})$ is invertible, and its inverse can be bounded by a constant independent of $T$ and $\varepsilon$.
\end{lemma}

\section{Proof of Theorem \ref{Main-asymptotic}: the single cavity case}\label{single-cavity}
\setcounter{equation}{0}

In this section, we show the asymptotic analysis of the solution to \eqref{bvp} as $\varepsilon \to 0$ for the single cavity case, i.e.  $D=\varepsilon B + z$. We recall that $f$ in \eqref{bvp} is a heat point source, namely, $f(x,\,t)=\Phi(x,\,t;\,z^*,\,0)$ with $z^*\in \mathbb R^3\setminus\overline D$. Express the solution of \eqref{bvp} as a single-layer heat potential
\begin{equation}\label{sol_bvp}
u(x,\,t) = \int_0^t \int_{\partial D} \Phi(x,\,t;\,y,\,\tau)\,\sigma(y,\,\tau)\,ds(y)d\tau,\quad (x,\,t)\in (\mathbb R^3\setminus\overline D)_T.
\end{equation}
In terms of the boundary condition in \eqref{bvp}, the density function $\sigma$ should satisfy
\begin{equation}\label{ie_bvp}
\mathcal S_{(\partial D)_T}[\sigma](x,\,t) = \Phi(x,\,t;\,z^*,\,0),\quad (x,\,t)\in (\partial D)_T.
\end{equation}
Since $\mathcal S_{(\partial D)_T}:\,L^2((\partial D)_T) \to H^{1,1/2}((\partial D)_T)$ is invertible, it follows that
\begin{equation*}
\sigma(x,\,t) = \mathcal S_{(\partial D)_T}^{-1}[\Phi_{(z^*,\,0)}](x,\,t),\quad (x,\,t)\in (\partial D)_T.
\end{equation*}
By Lemma \ref{S_ep} and \eqref{SD_3}, we derive that
\begin{eqnarray}\label{sigma_es1}
\|\sigma\|_{L^2((\partial D)_T)} & \leq &  \left \|\mathcal S_{(\partial D)_T}^{-1} \right\|_{\mathcal L \left(H^{1,1/2}((\partial D)_T),\,L^2((\partial D)_T)\right)} \|\Phi_{(z^*,\,0)} \|_{H^{1,1/2}((\partial D)_T)} \nonumber\\
& \leq & \varepsilon^{-1}
\left \| \mathcal S_{(\partial B)_{T_\varepsilon}}^{-1} \right\|_{\mathcal L \left(H^{1,1/2}((\partial B)_{T_\varepsilon}),\,L^2((\partial B)_{T_\varepsilon})\right)} \|\Phi_{(z^*,\,0)} \|_{H^{1,1/2}((\partial D)_T)} \nonumber\\
& \lesssim & \varepsilon^{-1} |\partial D|^{1/2} \leq C.
\end{eqnarray}
This also leads to
\begin{equation}\label{sigma_es2}
\int_0^t \int_{\partial D} |\sigma(y,\,\tau)|\,ds(y)d\tau \leq \|1\|_{L^2((\partial D)_T)}\cdot \|\sigma\|_{L^2((\partial D)_T)} \lesssim \varepsilon.
\end{equation}

In our argument, we also need to estimate $\partial_t \sigma$. Since $\Phi(x,\,t;\,z^*,\,0)$ is infinitely smooth with respect to $(x,\,t)\in (\partial D)_T$ and goes to zero identically as $t\to 0$, we derive from Lemma \ref{S_inv2} that $\sigma(\cdot,\,0)=0$ identically. Then, by taking the derivative of \eqref{sol_bvp} with respect to $t$ and using integration by parts, we get
\begin{eqnarray*}
\partial_t u(x,\,t) &=& \int_0^t \int_{\partial D} \partial_t \Phi(x,\,t;\,y,\,\tau)\,\sigma(y,\,\tau)\,ds(y)d\tau\\
&=& - \int_0^t \int_{\partial D} \partial_\tau \Phi(x,\,t;\,y,\,\tau)\,\sigma(y,\,\tau)\,ds(y)d\tau\\
&=& \int_0^t \int_{\partial D} \Phi(x,\,t;\,y,\,\tau)\,\partial_\tau\sigma(y,\,\tau)\,ds(y)d\tau, \quad (x,\,t)\in (\mathbb R^3\setminus\overline D)_T.
\end{eqnarray*}
Again, we obtain from the boundary condition in \eqref{bvp} that
\begin{equation}\label{sol_bvp_t}
\int_0^t \int_{\partial D} \Phi(x,\,t;\,y,\,\tau)\,\partial_\tau\sigma(y,\,\tau)\,ds(y)d\tau = \partial_t \Phi(x,\,t;\,z^*,\,0),\quad (x,\,t)\in (\partial D)_T.
\end{equation}
By the same derivations as for \eqref{sigma_es1} and \eqref{sigma_es2}, we have
\begin{equation}\label{sigma_es1_t}
\|\partial_t\sigma\|_{L^2((\partial D)_T)} \leq C
\end{equation}
and then
\begin{equation}\label{sigma_es2_t}
\int_0^t \int_{\partial D} |\partial_\tau\sigma(y,\,\tau)|\,ds(y)d\tau \lesssim \varepsilon.
\end{equation}

We now show the asymptotic expansion of the solution to \eqref{bvp} as $\varepsilon \to 0$, based on the integral representation \eqref{sol_bvp}. First, we observe that, for any fixed $(x,\,t)\in (\mathbb R^3\setminus\overline D)_T$, the function $\Phi(x,\,t;\,y,\,\tau)$ is sufficiently smooth with respect to $(y,\,\tau)\in\partial D \times (0,\,t)$. It follows from Taylor's expansion that
\begin{equation*}
|\Phi(x,\,t;\,y,\,\tau)-\Phi(x,\,t;\,z,\,\tau)|\lesssim \varepsilon
\end{equation*}
for $z\in D,\,y\in\partial D,\,x\in\mathbb R^3\setminus\overline D$. Then we obtain
\begin{equation*}
\left|\int_0^t \int_{\partial D} [\Phi(x,\,t;\,y,\,\tau)-\Phi(x,\,t;\,z,\,\tau)]\,\sigma(y,\,\tau)\,ds(y)d\tau \right| \lesssim \varepsilon^2,
\end{equation*}
and hence
\begin{eqnarray}\label{asym_1}
u(x,\,t) &=& \int_0^t \int_{\partial D} \Phi(x,\,t;\,z,\,\tau)\,\sigma(y,\,\tau)\,ds(y)d\tau \nonumber \\
& & + \int_0^t \int_{\partial D} [\Phi(x,\,t;\,y,\,\tau)-\Phi(x,\,t;\,z,\,\tau)]\,\sigma(y,\,\tau)\,ds(y)d\tau \nonumber \\
&=& \int_0^t \Phi(x,\,t;\,z,\,\tau) \left(\int_{\partial D} \sigma(y,\,\tau)\,ds(y)\right)d\tau + O(\varepsilon^2), \quad (x,\,t)\in (\mathbb R^3\setminus\overline D)_T.
\end{eqnarray}

\medskip
To derive the asymptotic expansion for the integral $\int_{\partial D} \sigma(y,\,t)\,ds(y)$, we rewrite the single-layer heat potential \eqref{sol_bvp} as
\begin{eqnarray}\label{ie_bvp2}
& &\int_0^t \int_{\partial D} \Phi(x,\,t;\,y,\,\tau) \, \sigma(y,\,\tau)\,ds(y)d\tau  \nonumber \\
&=& \int_{\partial D} \frac{1}{4\pi |x-y|} \left[ \int_0^t \frac{|x-y|}{2\sqrt{\pi}(t-\tau)^{3/2}}\,\exp\left( -\frac{|x-y|^2}{4(t-\tau)} \right)\, \sigma(y,\,\tau)\,d\tau \right] ds(y)
\end{eqnarray}
for $x\in\mathbb R^3\setminus\partial D$ and $t\in(0,\,T]$. Therefore, we can view the single-layer heat potential \eqref{sol_bvp} as a single-layer harmonic potential with the density
\begin{equation}\label{phi-x-y-t}
\varphi(x,\,y,\,t):= \int_0^t \frac{|x-y|}{2\sqrt{\pi}(t-\tau)^{3/2}}\,\exp\left( -\frac{|x-y|^2}{4(t-\tau)} \right)\, \sigma(y,\,\tau)\,d\tau, \quad x\in\mathbb R^3,\;y\in\partial D,\;t\in (0,\,T].
\end{equation}
By Lemma \ref{S_inv2}, the density $\sigma(x,\,t)$ in \eqref{ie_bvp} is continuous on $(\partial D)_T$, since $\Phi(x,\,t;\,z^*,\,0)$ is sufficiently smooth with respect to $(x,\,t)\in (\partial D)_T$. Based on this property, the continuity of the density function $\varphi$ is proved as follows.

\begin{lemma}\label{limit_phi}
If $\sigma(y,\,t)$ is continuous on $\partial D\times (0,\,T]$, the function $\varphi(x,\,y,\,t)$ is continuous on $\mathbb R^3\times \partial D\times (0,\,T]$ with
\begin{equation}\label{limit_1}
\lim_{x\to y}\varphi(x,\,y,\,t) = \sigma(y,\,t)
\end{equation}
for all $y\in\partial D$ and $t\in(0,\,T]$.
\end{lemma}

{\bf Proof.} Let
\begin{equation*}
\delta = \frac{|x-y|}{2\sqrt{t-\tau}}.
\end{equation*}
The direct calculations give
\begin{eqnarray}\label{def_I}
\varphi(x,\,y,\,t) &=& \int^{+\infty}_{\frac{|x-y|}{2\sqrt{t}}} \frac{2}{\sqrt{\pi}}\,\exp\left(-\delta^2\right)\,\sigma\left(y,\,t-\frac{|x-y|^2}{4\delta^2}\right)\,d\delta \nonumber\\
&=& \frac{2}{\sqrt{\pi}}\int_{\frac{|x-y|}{2\sqrt{t}}}^{\sqrt{|x-y|}} \,\exp\left(-\delta^2\right)\,\sigma\left(y,\,t-\frac{|x-y|^2}{4\delta^2}\right)\,d\delta  \nonumber\\
&& + \frac{2}{\sqrt{\pi}}\int_{\sqrt{|x-y|}}^{+\infty} \,\exp\left(-\delta^2\right)\,\left[\sigma\left(y,\,t-\frac{|x-y|^2}{4\delta^2}\right) - \sigma(y,\,t)\right]\,d\delta \nonumber\\
&& + \frac{2\sigma(y,\,t)}{\sqrt{\pi}}\int_{\sqrt{|x-y|}}^{+\infty} \,\exp(-\delta^2)\,d\delta \nonumber\\
&=:& I_1 + I_2 + I_3.
\end{eqnarray}

Clearly, we have
\begin{equation}\label{limit_2}
\lim_{x\to y} I_1(x,\,y,\,t) = 0
\end{equation}
uniformly on $\partial D$ and on compact subintervals of $(0,\,T]$. And also, we can easily see that
\begin{equation}\label{limit_3}
\lim_{x\to y} I_3(x,\,y,\,t) = \sigma(y,\,t) \frac{2}{\sqrt{\pi}} \int_0^{+\infty} \exp\left(-\delta^2\right)\,d\delta = \sigma(y,\,t)
\end{equation}
uniformly on $\partial D\times (0,\,T]$.

Let us consider $I_2(x,\,y,\,t)$. Since $\sigma(y,\,s)$ is continuous on $\partial D \times (0,\,T]$, for any $\epsilon>0$ there exists a positive constant $\delta_0$ such that $|\sigma(y,\,t_1) - \sigma(y,\,t_2)|\leq \epsilon$ for all $t_1$ and $t_2$ with $|t_1-t_2|<\delta_0$. Then, for all $|x-y|<4\delta_0$ and all $\delta \geq \sqrt{|x-y|}$, we have
\begin{equation*}
\frac{|x-y|^2}{4\delta^2} \leq \frac{|x-y|}{4} < \delta_0,
\end{equation*}
and hence
\begin{equation*}
\left| \sigma\left(y,\,t-\frac{|x-y|^2}{4\delta^2}\right) - \sigma(y,\,t) \right| \leq \epsilon.
\end{equation*}
Consequently, we obtain
\begin{equation*}
|I_2(x,\,y,\,t)| \leq \epsilon \frac{2}{\sqrt{\pi}} \int_0^{+\infty} \exp\left(-\delta^2\right)\,d\delta \leq \epsilon,
\end{equation*}
which implies that
\begin{equation}\label{limit_4}
\lim_{x\to y}I_2(x,\,y,\,t)=0.
\end{equation}
Combining \eqref{limit_2}, \eqref{limit_3} and \eqref{limit_4} yields \eqref{limit_1}. The proof is complete. \hfill $\Box$

Furthermore, we can estimate $|\sigma(y,\,t)-\varphi(x,\,y,\,t)|$ for $x,\,y\in \partial D$ and $t\in(0,\,T]$.
This property is crucial in the next steps.
\begin{lemma}\label{I}
Let $\sigma(y,\,t)$ be the solution to \eqref{ie_bvp}. Then we have
\begin{equation}\label{est_I}
\varphi(x,\,y,\,t)-\sigma(y,\,t)=  O \left(\vert x-y\vert \; \|\partial_t\sigma(y,\,\cdot)\|_{L^2(0,\,t)}\right)
\end{equation}
for $x,\, y$ such that $ \vert x-y\vert \ll 1$ and $t\in(0,\,T]$ uniformly with respect to $D$.
\end{lemma}

{\bf Proof.} We start from the formula
$$\varphi(x,\,y,\,t)=\int^{+\infty}_{\frac{\vert x-y\vert}{2\sqrt{t}}}\frac{2}{\sqrt{\pi}}\,\exp(-\delta^2)\,\sigma\left(y,\,t-\frac{|x-y|^2}{4\delta^2}\right)\,d\delta,$$
which we rewrite as
$$
\varphi(x,\,y,\,t)=\int^{+\infty}_{\frac{\vert x-y\vert}{2\sqrt{t}}}\frac{2}{\sqrt{\pi}}\,\exp(-\delta^2)\,\left[\sigma\left(y,\,t-\frac{|x-y|^2}{4\delta^2}\right)-\sigma(y,\,t)\right]\,d\delta+\sigma(y,\,t)\int^{+\infty}_{\frac{\vert x-y\vert}{2\sqrt{t}}}\frac{2}{\sqrt{\pi}}\exp(-\delta^2)\, d\delta.
$$
Thus, we obtain
\begin{eqnarray*}
\varphi(x,\,y,\,t)- \sigma(y,\,t) &=& \int^{+\infty}_{\frac{\vert x-y\vert}{2\sqrt{t}}}\frac{2}{\sqrt{\pi}}\,\exp(-\delta^2)\left[\sigma\left(y,\, t-\frac{|x-y|^2}{4\delta^2}\right)-\sigma(y,\, t)\right]\,d\delta \\
&& -\sigma(y,\,t)\int^{\frac{\vert x-y\vert}{2\sqrt{t}}}_{0}\frac{2}{\sqrt{\pi}}\,\exp(-\delta^2) \,d\delta.
\end{eqnarray*}
As $\sigma(y,\,0)=0$, then
$\sigma(y,\,t)=\int^{t}_0\partial_t \sigma(y,\,s)\,ds$ and by Cauchy-Schwartz inequality, we derive the estimate
$$\sigma(y,\,t)=O\left(t^{1/2}\,\Vert \sigma_t(y,\, \cdot) \Vert_{L^2(0,\,t)}\right).$$
Similarly, we have
$$\sigma\left(y,\,t-\frac{|x-y|^2}{4\delta^2}\right)-\sigma(y,\,t)=\int_t^{t-\frac{|x-y|^2}{4\delta^2}}\partial_t
\sigma(y,\,s)\,ds.$$
As $\delta \geq \frac{\vert x-y\vert}{2\sqrt{t}}$, then $t-\frac{\vert x-y\vert^2}{4\delta^2}\geq 0$, and hence
$$\int^{t-\frac{|x-y|^2}{4\delta^2}}_{t} \partial_t \sigma(y,\,s) \,ds = O\left(\frac{\vert x-y\vert^2}{4\delta^2}\,\Vert \partial_t \sigma(y,\,\cdot)\Vert_{L^2(0,\,t)}\right).$$
This means
\begin{eqnarray*}
\varphi(x,\,y,\,t)- \sigma(y,\,t) &=& O\left(\int^{+\infty}_{\frac{\vert x-y\vert}{2\sqrt{t}}}\frac{\exp(-\delta^2)}{\delta^2}\, d\delta\; \Vert \partial_t \sigma(y,\,\cdot)\Vert_{L^2(0,\,t)} \,\vert x-y\vert^2\right) \\
&& + O\left(\int^{\frac{\vert x-y\vert}{2\sqrt{t}}}_{0}\exp(-\delta^2)\, d\delta\, \Vert \partial_t \sigma(y,\,\cdot)\Vert_{L^2(0,\,t)}\, t^{1/2}\right).
\end{eqnarray*}
However, we have
$$\int^{\frac{\vert x-y\vert}{2\sqrt{t}}}_{0}\exp(-\delta^2) \,d\delta = O\left(\frac{\vert x-y\vert}{\sqrt{t}}\right),\quad \int^{+\infty}_{\frac{\vert x-y\vert}{2\sqrt{t}}}\frac{\exp(-\delta^2)}{\delta^2}\, d\delta = O\left(\frac{t}{|x-y|}\right).$$
Hence, we deduce
$$
\varphi(x,\,y,\,t)- \sigma(y,\,t) =O\left(\vert x-y\vert\,\Vert \partial_t \sigma(y,\,\cdot)\Vert_{L^2(0,\,t)} \right).
$$
The proof is complete. \hfill $\Box$

\bigskip
We are in a position to show the asymptotic expansion of the solution to \eqref{bvp} for the single cavity case.
\begin{theorem}\label{th_asym1}
For $x\in \mathbb R^3\setminus\overline D$ and $t\in(0,\,T]$, the solution $u(x,\,t)$ to \eqref{bvp} has the following asymptotic expansion:
\begin{equation}\label{asym_0}
u(x,\,t) = C_0\,\int_0^t \Phi(x,\,t;\,z,\,\tau)\, \Phi(z,\,\tau;\,z^*,\,0)\, d\tau + O(\varepsilon^2) \int_0^t \Phi(x,\,t;\,z,\,\tau) \, d\tau + O(\varepsilon^2), \quad \varepsilon\to 0
\end{equation}
with the constant $C_0$ defined by
\begin{equation*}
C_0:= \int_{\partial D} S_{\partial D}^{-1}[1](y)\,ds(y),
\end{equation*}
where $S_{\partial D}^{-1}$ is the inverse of the single-layer potential operator $S_{\partial D}$ corresponding to the Laplace equation, namely,
\begin{equation*}
S_{\partial D} [\varphi](x):= \int_{\partial D} \frac{1}{4\pi |x-y|}\, \varphi(y)\,ds(y).
\end{equation*}
\end{theorem}

{\bf Proof.} From \eqref{ie_bvp} and \eqref{ie_bvp2}, we obtain that
\begin{eqnarray}\label{asym_2}
\int_{\partial D} \frac{1}{4\pi |x-y|}\, \sigma(y,\,t) \,ds(y) &=& \Phi(x,\,t;\,z^*,\,0) + \int_{\partial D} \frac{1}{4\pi |x-y|}\,\left[ \sigma(y,\,t)-\varphi(x,\,y,\,t)\right]\,ds(y)
\end{eqnarray}
for $(x,\,t)\in (\partial D)_T$. Note that
\begin{equation*}
\Phi(x,\,t;\,z^*,\,0) = \Phi(z,\,t;\,z^*,\,0) + O(\varepsilon)
\end{equation*}
and, by Lemma \ref{I}, we have
\begin{equation*}
\int_{\partial D} \frac{1}{4\pi |x-y|}\,\left[ \sigma(y,\,t)-\varphi(x,\,y,\,t)\right]\,ds(y) = O\left( \varepsilon^{2}\right)
\end{equation*}
hold for $x\in\partial D$ and $t\in(0,\,T]$. Let $\sigma_z$ be the unique solution of
\begin{equation}\label{asym_3}
\int_{\partial D} \frac{1}{4\pi |x-y|}\, \sigma_z(y,\,t) \,ds(y) = \Phi(z,\,t;\,z^*,\,0).
\end{equation}
Then we have
\begin{equation}\label{asym_4}
\int_{\partial D} \frac{1}{4\pi |x-y|}\, \left[\sigma(y,\,t) - \sigma_z(y,\,t)\right]\,ds(y) = O(\varepsilon) + O\left(\varepsilon^{2}\right) = O(\varepsilon), \mbox{ uniformly for } x \in \partial D.
\end{equation}
We recall that $S_{\partial D}$ from $H^{-1}(\partial D)$ to $ L^2(\partial D)$ is an isomorphism and $\Vert S_{\partial D}^{-1}\Vert_{\mathcal L(L^2(\partial D), H^{-1}(\partial D))}=O(\varepsilon^{-1})$.
This implies that
\begin{equation}\label{asym_5}
\| \sigma - \sigma_z \|_{H^{-1}(\partial D)} = \left\|S_{\partial D}^{-1}\right\|_{\mathcal L\left(L^2(\partial D),\,H^{-1}(\partial D) \right)} O(\varepsilon^2) = O(\varepsilon),
\end{equation}
and then
\begin{equation}\label{asym_6}
\int_{\partial D} \sigma(y,\,t)\,ds(y) = \int_{\partial D} \sigma_z(y,\,t)\,ds(y) + O(\varepsilon^2).
\end{equation}

From \eqref{asym_3}, we know that
\begin{equation*}
\sigma_z(x,\,t) = S_{\partial D}^{-1}[1](x)\, \Phi(z,\,t;\,z^*,\,0),
\end{equation*}
and therefore
\begin{equation}\label{asym_7}
\int_{\partial D}\sigma_z(y,\,t)\,ds(y) = C_0\,\Phi(z,\,t;\,z^*,\,0).
\end{equation}
That is,
\begin{equation}\label{asym_7_2}
\int_{\partial D}\sigma(y,\,t)\,ds(y) = C_0\,\Phi(z,\,t;\,z^*,\,0) + O(\varepsilon^2).
\end{equation}
Then, by inserting it into \eqref{asym_1}, we have, for $x\in \mathbb R^3\setminus\overline D$ and $t\in (0, \;T]$, that
\begin{equation}\label{asym_8}
u(x,\,t) = C_0\,\int_0^t \Phi(x,\,t;\,z,\,\tau)\, \Phi(z,\,\tau;\,z^*,\,0)\, d\tau + O(\varepsilon^2) \int_0^t \Phi(x,\,t;\,z,\,\tau) \, d\tau + O(\varepsilon^2).
\end{equation}
The proof is now complete. \hfill $\Box$

\begin{remark}
Observe that when the source of the heat $z^*$ and the location of the receiver $x$ are away from the cavity $D$ then for $t \sim \varepsilon$, the dominating term $C_0\,\int_0^t \Phi(x,\,t;\,z,\,\tau)\, \Phi(z,\,\tau;\,z^*, \,0)\, d\tau$ behaves as $\varepsilon^2$ and hence it is lost in the error term. But this is not a surprise. However, when $z^*$ or/and $x$ are close to the cavity $D$, then the first term stays a dominating term even for very short time $t$, i.e. $t \sim \varepsilon$.
\end{remark}

\section{Proof of Theorem \ref{Main-asymptotic}: the multiple cavities case}\label{multiple-cavities}

In this section, we show the asymptotic analysis of the solution to \eqref{bvp} as $\varepsilon \to 0$ for the multiple cavities case. We recall that $D_j := \varepsilon B_j + z_j$ to be small cavities characterized by the parameter $\varepsilon >0$ and the locations $z_j\in\mathbb R^3,\,j=1,\,2,\,\cdots,\,M$.

\subsection{Integral representation of the solution}
\setcounter{equation}{0}

We express the solution to \eqref{bvp} as a single-layer heat potential
\begin{equation}\label{sol_bvp_M}
u(x,\,t) = \sum_{j=1}^M\int_0^t \int_{\partial D_j} \Phi(x,\,t;\,y,\,\tau)\,\sigma_j(y,\,\tau)\,ds(y)d\tau,\quad (x,\,t)\in (\mathbb R^3\setminus\overline D)_T,
\end{equation}
where $\sigma_j,\,j=1,\,2,\,\cdots,\,M$ are density functions to be determined. In terms of the boundary condition, the density functions should satisfy
\begin{equation}\label{ie_bvp_M}
\sum_{j=1}^M\mathcal S_{(\partial D_j)_T}[\sigma_j](x,\,t) = \Phi(x,\,t;\,z^*,\,0),\quad (x,\,t)\in (\partial D)_T.
\end{equation}
Define the operator $\mathcal S_{ji}$ by
\begin{equation}\label{S_ij}
\mathcal S_{ji}[\sigma_j](x,\,t) := \int_0^t\int_{\partial D_j}\Phi(x,\,t;\,y,\,\tau)\,\sigma_j(y,\,\tau)\,ds(y)d\tau, \quad (x,\,t)\in (\partial D_i)_T.
\end{equation}
Then \eqref{ie_bvp_M} can be rewritten as
\begin{equation}\label{ie_bvp_M_1}
\mathcal S_{ii}[\sigma_i] + \sum_{\substack{j=1\\j\neq i}}^M \mathcal S_{ji}[\sigma_j] = \Phi(x,\,t;\,z^*,\,0),\quad (x,\,t)\in (\partial D_i)_T,\quad i=1,\,2,\,\cdots,\,M
\end{equation}
or
\begin{equation}\label{ie_bvp_M_2}
\sigma_i + \mathcal S_{ii}^{-1} \sum_{\substack{j=1\\j\neq i}}^M \mathcal S_{ji}[\sigma_j] = \mathcal S_{ii}^{-1} \left[ \Phi_{(z^*,\,0)}\right],\quad (x,\,t)\in (\partial D_i)_T,\quad i=1,\,2,\,\cdots,\,M.
\end{equation}

The unique solvability of the system \eqref{ie_bvp_M_1}, or \eqref{ie_bvp_M_2}, can be justified using standard Fredholm alternative. Indeed, the operators $\mathcal S_{ii}: L^2((\partial D_i)_T) \to H^{1,1/2}((\partial D_i)_T)$ are invertible, while $\mathcal S_{ji}, j\neq i,$ are compact ones. To show uniqueness of the solution of problem \eqref{ie_bvp_M_1}, we use uniqueness of the problem (\ref{bvp}) and the jump relations of the adjoint of the double layer operator. We omit the details here.

The next step is to derive a priori estimates of the densities $\sigma_i$'s.

\subsection{A priori estimate of the densities $\sigma_i$'s}

We start with the following singularity properties related to the heat fundamental solution.

\begin{lemma}\label{int_phi}
For $x\not=y$ and $|x-y|\to 0$, we have
\begin{eqnarray}
&&\left(\int_0^T\int_0^T |\Phi(x,\,t;\,y,\,\tau)|^2\,d\tau dt\right)^{1/2}=O(|x-y|^{-2}), \label{int_phi_e1}\\
&&\left(\int_0^T\int_0^T |\partial_{x_i}\Phi(x,\,t;\,y,\,\tau)|^2\,d\tau dt\right)^{1/2}=O(|x-y|^{-3}), \quad i=1,\,2,\,3,\label{int_phi_e2}\\
&&\left(\int_0^T\int_0^T |D_t^{1/2}\Phi(x,\,t;\,y,\,\tau)|^2\,d\tau dt\right)^{1/2}=O(|x-y|^{-3}). \label{int_phi_e3}
\end{eqnarray}
\end{lemma}

{\bf Proof.} Let $\zeta=|x-y|$ and $\delta=\frac{\zeta}{2\sqrt{t-\tau}}$. By direct calculations, we have
\begin{eqnarray}\label{int_phi_1}
\int_0^T\int_0^T |\Phi(x,\,t;\,y,\,\tau)|^2\,d\tau dt
&=& \int_0^T\int_0^t \frac{1}{[4\pi (t-\tau)]^3} \exp\left(- \frac{|x-y|^2}{2(t-\tau)} \right)\,d\tau dt\nonumber\\
&=& \frac{\zeta^{-4}}{2\pi^3}\int_0^T dt \int_{\frac{\zeta}{2\sqrt{t}}}^{+\infty} \delta^3 \exp(-2\delta^2)\,d\delta\nonumber\\
&=&  \frac{\zeta^{-4}}{16\pi^3} \int_0^T \left[ \frac{\zeta^2}{2t}\exp(-\frac{\zeta^2}{2t}) + \exp(-\frac{\zeta^2}{2t})\right]dt.
\end{eqnarray}
Since
\begin{eqnarray}\label{int_phi_2}
\int_0^T \exp(-\frac{\zeta^2}{2t})\,dt = \left[ t \exp(-\frac{\zeta^2}{2t})\right]_0^T - \int_0^T \frac{\zeta^2}{2t}\exp(-\frac{\zeta^2}{2t})\,dt = T \exp(-\frac{\zeta^2}{2T}) - \int_0^T \frac{\zeta^2}{2t}\exp(-\frac{\zeta^2}{2t})\,dt,
\end{eqnarray}
we obtain from \eqref{int_phi_1} that
\begin{eqnarray}\label{int_phi_3}
\int_0^T\int_0^T |\Phi(x,\,t;\,y,\,\tau)|^2\,dtd\tau = \frac{T\zeta^{-4}}{16\pi^3} \exp(-\frac{\zeta^2}{2T}) \sim \frac{T\zeta^{-4}}{16\pi^3} \quad \textrm{as}\;\zeta\to 0.
\end{eqnarray}
This completes the proof of \eqref{int_phi_e1}.

\medskip
Next, let us prove \eqref{int_phi_e2}. The direct calculations give
\begin{eqnarray}\label{int_phi_4}
\int_0^T\int_0^T |\partial_{x_i}\Phi(x,\,t;\,y,\,\tau)|^2\,d\tau dt
&\leq& \int_0^T\int_0^t \frac{|x-y|^2}{[4\pi (t-\tau)]^5} \exp\left(- \frac{|x-y|^2}{2(t-\tau)} \right)\,d\tau dt\nonumber\\
&=& \frac{\zeta^{-6}}{2\pi^5}\int_0^T dt \int_{\frac{\zeta}{2\sqrt{t}}}^{+\infty} \delta^7 \exp(-2\delta^2)\,d\delta.
\end{eqnarray}
By repeatedly using integration by parts, we have
\begin{equation}\label{int_phi_5}
\int_{\frac{\zeta}{2\sqrt{t}}}^{+\infty} \delta^7 \exp(-2\delta^2)\,d\delta = \frac{\zeta^6}{2^8t^3}\exp(-\frac{\zeta^2}{2t}) + \frac{3\zeta^4}{2^7t^2}\exp(-\frac{\zeta^2}{2t}) + \frac{3}{16}\left[ \frac{\zeta^2}{2t}\exp(-\frac{\zeta^2}{2t}) + \exp(-\frac{\zeta^2}{2t})\right].
\end{equation}
It follows that
\begin{equation}\label{int_phi_6}
\int_0^T dt \int_{\frac{\zeta}{2\sqrt{t}}}^{+\infty} \delta^7 \exp(-2\delta^2)\,d\delta=O(1), \quad \zeta\to 0.
\end{equation}
Then the estimate \eqref{int_phi_e2} comes from \eqref{int_phi_4} and \eqref{int_phi_6}.

\medskip
Finally, let us show \eqref{int_phi_e3}. Due to the estimate
\begin{equation}\label{int_phi_7}
|D_t^{1/2}\Phi(x,\,t;\,y,\,\tau)| \lesssim \frac{1}{(t-\tau)^{3/2}} \,\frac{1}{|x-y|} \,\exp\left( -c_0\frac{|x-y|^2}{t-\tau} \right)
\end{equation}
with a positive constant $c_0$; see \cite[Lemma A.1]{Brown1987}, we get
\begin{eqnarray}\label{int_phi_8}
\int_0^T\int_0^T |D_t^{1/2}\Phi(x,\,t;\,y,\,\tau)|^2\,d\tau dt
&\lesssim & \int_0^T\int_0^t \frac{1}{(t-\tau)^3}\,\frac{1}{|x-y|^2} \exp\left(- 2c_0\frac{|x-y|^2}{t-\tau} \right)\,d\tau dt\nonumber\\
&=& 2^5\zeta^{-6}\int_0^T dt \int_{\frac{\zeta}{2\sqrt{t}}}^{+\infty} \delta^3 \exp(-8c_0\delta^2)\,d\delta \nonumber\\
&=& O(\zeta^{-6}).
\end{eqnarray}
The proof is now complete. \hfill $\Box$

\begin{lemma}\label{int_dphi}
For $x\not=y$ and $|x-y|\to 0$, we have
\begin{equation}\label{int_dphi_e}
\left(\int_0^t \left|\nabla_x\Phi(x,\,t;\,y,\,\tau)\right|^2\,d\tau\right)^{1/2} = O(|x-y|^{-3}), \quad t\in(0,\,T].
\end{equation}
\end{lemma}

{\bf Proof.} Let $\zeta=|x-y|$ and $\delta=\frac{\zeta}{2\sqrt{t-\tau}}$. We have
\begin{equation}\label{int_dphi_1}
\int_0^t \left|\nabla_x\Phi(x,\,t;\,y,\,\tau)\right|^2\,d\tau \leq \int_0^t \frac{\zeta^2}{[4\pi(t-\tau)]^5}\exp\left(-\frac{\zeta^2}{2(t-\tau)}\right)d\tau = \frac{\zeta^{-6}}{2\pi^5}\int_{\frac{\zeta}{2\sqrt{t}}}^{+\infty} \delta^7 \exp(-2\delta^2)\,d\delta,
\end{equation}
which gives \eqref{int_dphi_e} by using \eqref{int_phi_5}. The proof is complete. \hfill $\Box$

\bigskip
To proceed, we recall the single-layer operator $S_{\partial D_i}$:
\begin{equation}\label{S_Di}
S_{\partial D_i}[\sigma_i](x,\,t) := \int_{\partial D_i} \frac{1}{4\pi |x-y|}\,\sigma_i(y,\,t)\,ds(y), \quad (x,\,t)\in (\partial D_i)_T,
\end{equation}
and write
\begin{equation*}\label{MS_1}
\mathcal S_{ii}[\sigma_i] = S_{\partial D_i}[\sigma_i] + S_{\partial D_i}[\varphi_i - \sigma_i],
\end{equation*}
where
\begin{equation*}\label{phi_i}
\varphi_i(x,\,y,\,t):= \int_0^t \frac{|x-y|}{2\sqrt{\pi}(t-\tau)^{3/2}}\exp\left(-\frac{|x-y|^2}{4(t-\tau)}\right)\sigma_i(y,\,\tau)\,d\tau, \quad x\in\mathbb R^3,\,y\in\partial D_i,\,t\in(0,\,T].
\end{equation*}
In addition, we define
\begin{equation}\label{def_qi}
q_j(t) :=  \int_{\partial D_j} \sigma_j(y,\,t)\,ds(y),
\end{equation}
and write
\begin{equation*}
\mathcal S_{ji}[\sigma_j](x,\,t) = \int_0^t \Phi(z_i,\,t;\,z_j,\,\tau)\,q_j(\tau)\,d\tau + A_{ji} + B_{ji}
\end{equation*}
with
\begin{eqnarray*}
A_{ji} &:=& \int_0^t \left[ \Phi(x,\,t;\,z_j,\,\tau) - \Phi(z_i,\,t;\,z_j,\,\tau) \right] q_j(\tau)\,d\tau,\\
B_{ji} &:=& \int_0^t\int_{\partial D_j} \left[ \Phi(x,\,t;\,y,\,\tau) - \Phi(x,\,t;\,z_j,\,\tau) \right] \sigma_j(y,\,\tau)\,ds(y)d\tau.
\end{eqnarray*}
Hence, \eqref{ie_bvp_M_1} becomes
\begin{equation}\label{MS2}
S_{\partial D_i}[\sigma_i] =- \sum_{\substack{j=1\\j\neq i}}^M \int_0^t \Phi(z_i,\,t;\,z_j,\,\tau)\,q_j(\tau)\,d\tau +\Phi_{(z^*,\,0)} - S_{\partial D_i}[\varphi_i - \sigma_i] - \sum_{\substack{j=1\\j\neq i}}^M \left(A_{ji} -  B_{ji}\right).
\end{equation}

\medskip
To deduce the estimate of $\sigma_i$'s, we decompose $A_{ji}$ as
\begin{equation*}\label{MS_A}
A_{ji} = A_{ji}^{(1)} + A_{ji}^{(2)}
\end{equation*}
with
\begin{equation}\label{MS_A1}
A_{ji}^{(1)} := \int_0^t (x-z_i)\cdot \nabla_x \Phi(z_i,\,t;\,z_j,\,\tau)\, q_j(\tau)\,d\tau
\end{equation}
and
\begin{equation}\label{MS_A2}
A_{ji}^{(2)} := \frac{1}{2}\sum_{k,\,l=1}^3 \int_0^t (x_k-{z_i}_k)\,(x_l-{z_i}_l)\,\partial_{x_k}\partial_{x_l} \Phi(z_i^*,\,t;\,z_j,\,\tau)\, q_j(\tau)\,d\tau, \; z_i^*=z_i+\theta(x-z_i),\,0<\theta<1.
\end{equation}
Then we decompose $\sigma_i$ as $\sigma_i=\sigma_i^{(1)}+\sigma_i^{(2)}+\sigma_i^{(3)}$, where $\sigma_i^{(1)},\,\sigma_i^{(2)}$ and $\sigma_i^{(3)}$ are defined by
\begin{equation}\label{MS2_1}
S_{\partial D_i}[\sigma_i^{(1)}] = - \sum_{\substack{j=1\\j\neq i}}^M A_{ji}^{(1)},
\end{equation}
\begin{equation}\label{MS2_2}
S_{\partial D_i}[\sigma_i^{(2)}] = - S_{\partial D_i}[\varphi_i - \sigma_i],
\end{equation}
and
\begin{equation}\label{MS2_3}
S_{\partial D_i}[\sigma_i^{(3)}] = - \sum_{\substack{j=1\\j\neq i}}^M \int_0^t \Phi(z_i,\,t;\,z_j,\,\tau)\,q_j(\tau)\,d\tau +\Phi_{(z^*,\,0)}  - \sum_{\substack{j=1\\j\neq i}}^M \left(A_{ji}^{(2)} -  B_{ji}\right),
\end{equation}
respectively.

First, we estimate $\sigma_i^{(1)}$ by
\begin{eqnarray}\label{MS2_11}
\|\sigma_i^{(1)} \|_{L^2(\partial D_i)} &\lesssim& \sum_{\substack{j=1\\j\neq i}}^M \|S_{\partial D_i}^{-1}[x-z_i] \|_{L^2(\partial D_i)}\,\int_0^t |\nabla_x \Phi(z_i,\,t;\,z_j,\,\tau)|\,|q_j(\tau)|\,d\tau \nonumber\\
&\lesssim& \sum_{\substack{j=1\\j\neq i}}^M \varepsilon \left( \int_0^t |\nabla_x \Phi(z_i,\,t;\,z_j,\,\tau)|^2\,d\tau\right)^{1/2}\, \left( \int_0^t |q_j(\tau)|^2\,d\tau\right)^{1/2} \nonumber\\
&\lesssim& \sum_{\substack{j=1\\j\neq i}}^M \varepsilon^2\, d_{ji}^{-3}\, \| \sigma_j \|_{L^2((\partial D_j)_T)}.
\end{eqnarray}

Second, let us estimate $\sigma_i^{(2)}$. To this end, we need the following result:
\begin{lemma}\label{I_i}
Let $\sigma_i(y,\,t),\,i=1,\,2,\,\cdots,\,M$ be the solution to \eqref{ie_bvp_M}. Then we have
\begin{equation}\label{est_Ii}
\varphi_i(x,\,y,\,t)-\sigma_i(y,\,t)=  \sum_{n=1}^{+\infty} \alpha_n(t)\, \vert x-y\vert^n \,\sigma_i(y,\,t) + \sum_{n=1}^{+\infty} \left(\sum_{m=n}^{+\infty} \beta_{m,n}(t)\, |x-y|^m \right)\, \partial_t^n\sigma_i(y,\,t)
\end{equation}
for $x,\, y$ such that $ \vert x-y\vert \ll 1$ and $t\in(0,\,T]$ uniformly with respect to $D_i$, where $\alpha_n(t)$ and $\beta_{m,n}(t)$ are smooth in $(0,\,T]$.
\end{lemma}

{\bf Proof.} Recall that
\begin{eqnarray}\label{MS2_21}
\varphi_i(x,\,y,\,t)- \sigma_i(y,\,t) &=& \int^{+\infty}_{\frac{\vert x-y\vert}{2\sqrt{t}}}\frac{2}{\sqrt{\pi}}\,\exp(-\delta^2)\left[\sigma_i\left(y,\, t-\frac{|x-y|^2}{4\delta^2}\right)-\sigma_i(y,\, t)\right]\,d\delta \nonumber\\
&& -\sigma_i(y,\,t)\int^{\frac{\vert x-y\vert}{2\sqrt{t}}}_{0}\frac{2}{\sqrt{\pi}}\,\exp(-\delta^2) \,d\delta.
\end{eqnarray}
Note that
\begin{equation}\label{MS2_22}
\int^{\frac{\vert x-y\vert}{2\sqrt{t}}}_{0}\frac{2}{\sqrt{\pi}}\,\exp(-\delta^2) \,d\delta = \sum_{n=1}^{+\infty} \alpha_n \,\left( \frac{\vert x-y\vert}{2\sqrt{t}} \right)^n =:\sum_{n=1}^{+\infty} \alpha_n(t)\, \vert x-y\vert^n.
\end{equation}
In addition, for fixed $y$, we have
\begin{eqnarray*}
\sigma_i\left(y,\, t-\frac{|x-y|^2}{4\delta^2}\right)-\sigma_i(y,\, t) = \sum_{n=1}^{+\infty} \partial_t^n\sigma_i(y,\,t)\, \left(-\frac{|x-y|^2}{4\delta^2} \right)^{n},
\end{eqnarray*}
and then
\begin{eqnarray}\label{MS2_23}
&&\int^{+\infty}_{\frac{\vert x-y\vert}{2\sqrt{t}}}\frac{2}{\sqrt{\pi}}\,\exp(-\delta^2)\left[\sigma_i\left(y,\, t-\frac{|x-y|^2}{4\delta^2}\right)-\sigma_i(y,\, t)\right]\,d\delta \nonumber\\
&=& \sum_{n=1}^{+\infty} (-|x-y|^2)^n\, \partial_t^n\sigma_i(y,\,t)\,\int^{+\infty}_{\frac{\vert x-y\vert}{2\sqrt{t}}}\frac{2}{\sqrt{\pi}}\,\frac{\exp(-\delta^2)}{(4\delta^2)^n}\,d\delta \nonumber\\
&=& \sum_{n=1}^{+\infty} \left(\sum_{m=n}^{+\infty} \beta_{m,n}(t)\, |x-y|^m \right)\, \partial_t^n\sigma_i(y,\,t).
\end{eqnarray}
The proof is completed by inserting \eqref{MS2_22} and \eqref{MS2_23} into \eqref{MS2_21}. \hfill $\Box$

In view of Lemma \ref{I_i}, let us consider the following two equations:
\begin{eqnarray}
S_{\partial D_i}[\sigma_{n}^{(2)}] &=& S_{\partial D_i}[|x-y|^n\, \sigma_i(y,\,t)], \label{MS2_24}\\
S_{\partial D_i}[\sigma_{m,n}^{(2)}] &=& S_{\partial D_i}[|x-y|^m\,\partial_t^n \sigma_i(y,\,t)]. \label{MS2_25}
\end{eqnarray}
By scaling on the space variable, we get from \eqref{MS2_24} that
\begin{equation*}
\varepsilon\, S_{\partial B_i} [\tilde\sigma_{n}^{(2)}] = \varepsilon^{n+1}\,S_{\partial B_i}[|\xi-\eta|^2\,\tilde \sigma_i(\eta,\,t)],
\end{equation*}
where $x=\varepsilon\, \xi,\,y=\varepsilon\, \eta,\,\tilde \sigma_{n}^{(2)}(\eta,\,t)=\sigma_{n}^{(2)}(\varepsilon\eta,\,t)$ and $\tilde \sigma_i (\eta,\,t)=\sigma_i(\varepsilon\,\eta,\,t)$. Then we have
\begin{equation*}
\| \tilde \sigma_{n}^{(2)}\|_{L^2(\partial B_i)} \lesssim \varepsilon^n \|\tilde \sigma \|_{L^2(\partial B_i)},
\end{equation*}
and also
\begin{equation}\label{MS2_26}
\| \sigma_{n}^{(2)}\|_{L^2(\partial D_i)} \lesssim \varepsilon^n \| \sigma_i \|_{L^2(\partial D_i)}.
\end{equation}
Similarly, we obtain from \eqref{MS2_25} that
\begin{equation}\label{MS2_27}
\| \sigma_{n,n}^{(2)}\|_{L^2(\partial D_i)} \lesssim \varepsilon^n \| \partial_t^n\sigma_i \|_{L^2(\partial D_i)}.
\end{equation}
Hence, we have
\begin{eqnarray}\label{MS2_28}
\|\sigma_i^{(2)} \|_{L^2(\partial D_i)} &\lesssim& \sum_{n=1}^{+\infty} \varepsilon^n\,\left[\| \sigma_i \|_{L^2(\partial D_i)} + \| \partial_t^n\sigma_i \|_{L^2(\partial D_i)}\right] \nonumber\\
&=& O(\varepsilon\| \sigma_i \|_{L^2(\partial D_i)})+\sum_{n=1}^{+\infty} \varepsilon^n\,\| \partial_t^n\sigma_i \|_{L^2(\partial D_i)}.
\end{eqnarray}

Finally, let us estimate $\sigma_i^{(3)}$. Using the estimate of the gradient of the heat kernel; see (\ref{fun_est2}), we see that
if $z^*$ is away from $\cup_{i=1}^MD_i$, then
\begin{equation}
\Vert \Phi_{(z^*,\,0)}\Vert_{H^1(\partial D_i)}=O(\varepsilon).
\end{equation}
To proceed, we show the following estimates:

\begin{lemma}\label{lem_AB}
For $t\in(0,\,T]$ and $i,\,j=1,\,2,\,\cdots,\, M$ with $i\neq j$, we have
\begin{eqnarray}
\|A_{ji}^{(1)} \|^2_{H^1(\partial D_i)} \,  &\lesssim & \left( \varepsilon^6 d_{ji}^{-6} + \varepsilon^4 d_{ji}^{-6} \right) \|\sigma_j \|_{L^2((\partial D_j)_T)}^2, \label{est_A1}\\
\|A_{ji}^{(2)}\|^2_{H^1(\partial D_i)} \, &\lesssim & \left( \varepsilon^8 d_{ji}^{-8} + \varepsilon^6 d_{ji}^{-8} \right) \|\sigma_j \|_{L^2((\partial D_j)_T)}^2, \label{est_A2}\\
\|B_{ji}\|^2_{H^1(\partial D_i)}\, &\lesssim & \left( \varepsilon^6 d_{ji}^{-6} + \varepsilon^6 d_{ji}^{-8} \right) \|\sigma_j \|_{L^2((\partial D_j)_T)}^2. \label{est_B}
\end{eqnarray}
\end{lemma}

{\bf Proof.} We only prove the estimate \eqref{est_A1}, since the others can be shown in the same way. By direct calculations and Lemma \ref{int_dphi}, we obtain
\begin{eqnarray*}
\int_{\partial D_i} \left| A_{ji}^{(1)}\right|^2\,ds(x)
&=& \int_{\partial D_i}\left| \int_0^t (x-z_i)\cdot \nabla_x \Phi(z_i,\,t;\,z_j,\,\tau)\, q_j(\tau)\,d\tau \right|^2\,ds(x) \\
&\lesssim& \|\sigma_j \|_{L^2((\partial D_j)_T)}^2 \int_{\partial D_i}\int_{\partial D_j} \left( \int_0^t |x-z_i|^2\, |\nabla_x \Phi(z_i,\,t;\,z_j,\,\tau) |^2\,d\tau  \right)\,ds(x)ds(y)\\
&\lesssim& \varepsilon^6 d_{ji}^{-6} \|\sigma_j \|_{L^2((\partial D_j)_T)}^2
\end{eqnarray*}
and
\begin{eqnarray*}
\int_{\partial D_i} \left|\nabla_{\mathrm{tan}} A_{ji}^{(1)}\right|^2\,ds(x)
&\lesssim& \int_{\partial D_i} \left(\int_0^t \left|\nabla_x \Phi(z_i,\,t;\,z_j,\,\tau)\, q_j(\tau)\right|^2\,d\tau \right)\,ds(x)\\
&\lesssim& \|\sigma_j \|_{L^2((\partial D_j)_T)}^2 \int_{\partial D_i}\int_{\partial D_j} \left(\int_0^t |\nabla_x \Phi(z_i,\,t;\,z_j,\,\tau) |^2 \,d\tau \right)\,ds(x)ds(y)\\
&\lesssim& \varepsilon^4 d_{ji}^{-6} \|\sigma_j \|_{L^2((\partial D_j)_T)}^2.
\end{eqnarray*}
So, the estimate \eqref{est_A1} is justified.  \hfill $\Box$

\medskip
Note that the term $ \sum_{j\neq i} \int_0^t \Phi(z_i,\,t;\,z_j,\,\tau)\,q_j(\tau)\,d\tau $ is independent of the space variable, hence its $H^1(\partial D_i)$-norm reduces to its
$L^2(\partial D_i)$-norm. Then
\begin{eqnarray*}
\int_{\partial D_i} \left| \int_0^t \Phi(z_i,\,t;\,z_j,\,\tau)\,q_j(\tau)\,d\tau \right|^2\,ds(x)
&\lesssim&  \varepsilon^2\|q_j\|_{L^2}^2  \int_0^t \left| \Phi(z_i,\,t;\,z_j,\,\tau)\right|^2 \,d\tau\\
&\lesssim& \varepsilon^2 d_{ji}^{-4}\,\|q_j\|_{L^2}^2
\lesssim \varepsilon^4 d_{ji}^{-4} \|\sigma_j \|_{L^2((\partial D_j)_T)}^2.
\end{eqnarray*}

\bigskip
Then, by \eqref{MS2_3}, we get
\begin{equation*}\label{est_sigma-}
\|S_{\partial D_i} [\sigma_i^{(3)}]\|_{H^1(\partial D_i)}\lesssim \sum_{\substack{j=1\\j\neq i}}^M \varepsilon^2 d_{ji}^{-2} \|\sigma_j\|_{L^2((\partial D_j)_T)} + O(\varepsilon) + \sum_{\substack{j=1\\j\neq i}}^M  O\left(\varepsilon^3\, d_{ji}^{-4}\,\|\sigma_j \|_{L^2((\partial D_j)_T)} \right).
\end{equation*}
Using the property $\Vert S_{\partial D_i}^{-1} \Vert_{\mathcal L(H^1(\partial D_i),\,L^2(\partial D_i))} = O(\varepsilon^{-1})$ and integrating over $t$, we deduce that
\begin{equation}\label{est_sigma}
\| \sigma_i^{(3)}\|_{L^2((\partial D_i)_T)}
\lesssim \sum_{\substack{j=1\\j\neq i}}^M \varepsilon d_{ji}^{-2} \|\sigma_j\|_{L^2((\partial D_j)_T)} + O(1) + \sum_{\substack{j=1\\j\neq i}}^M  O\left(\varepsilon^2\, d_{ji}^{-4}\,\|\sigma_j \|_{L^2((\partial D_j)_T)} \right).
\end{equation}
Combing \eqref{MS2_11}, \eqref{MS2_28} and \eqref{est_sigma}, we have
\begin{eqnarray}\label{est_sigma_0}
\| \sigma_i\|_{L^2((\partial D_i)_T)} &\lesssim& \sum_{\substack{j=1\\j\neq i}}^M \varepsilon d_{ji}^{-2} \|\sigma_j\|_{L^2((\partial D_j)_T)} + O(1) + \sum_{\substack{j=1\\j\neq i}}^M  O\left(\varepsilon^2\, d_{ji}^{-4}\,\|\sigma_j \|_{L^2((\partial D_j)_T)} \right) \nonumber\\
&& + O(\varepsilon\| \sigma_i \|_{L^2(\partial D_i)})+\sum_{n=1}^{+\infty} \varepsilon^n\,\| \partial_t^n\sigma_i \|_{L^2(\partial D_i)}.
\end{eqnarray}

To estimate $\|\partial_t^n \sigma_i\|_{L^2((\partial D_i)_T)}$, we perform the time derivative for \eqref{ie_bvp_M_1}, and use integration by parts for the first term and the fact that $\sigma_i(\cdot,\,0)=0$. Then we get
\begin{equation}\label{sigmat_1}
\mathcal S_{ii}[\partial_t\sigma_i] + \sum_{\substack{j=1\\j\neq i}}^M \partial_t \left(\mathcal S_{ji}[\sigma_j]\right) = \partial_t \Phi(x,\,t;\,z^*,\,0),\quad (x,\,t)\in (\partial D_i)_T,
\end{equation}
or equivalently,
\begin{equation}\label{sigmat_2}
\partial_t\sigma_i = \mathcal S_{ii}^{-1}\Bigg[ \partial_t \Phi(x,\,t;\,z^*,\,0) - \sum_{\substack{j=1\\j\neq i}}^M \partial_t \left(\mathcal S_{ji}[\sigma_j]\right) \Bigg],\quad (x,\,t)\in (\partial D_i)_T.
\end{equation}
Since we have
\begin{equation}\label{sigmat_3}
\left\| \mathcal S_{ii}^{-1} \right\|_{\mathcal L\left(H^{1,1/2}((\partial D_i)_T),\,L^2((\partial D_i)_T)\right)} \lesssim \varepsilon^{-1},
\end{equation}
it suffices to estimate
\begin{equation*}
\Big\|\partial_t \Phi(x,\,t;\,z^*,\,0) - \sum_{\substack{j=1\\j\neq i}}^M \partial_t \left(\mathcal S_{ji}[\sigma_j]\right)\Big\|_{H^{1,1/2}((\partial D_i)_T)}.
\end{equation*}

First, it can be easily seen that
\begin{equation}\label{sigmat_4}
\Big\|\partial_t \Phi(x,\,t;\,z^*,\,0)\Big\|_{H^{1,1/2}((\partial D_i)_T)} \lesssim \varepsilon.
\end{equation}
Second, let us define $\mathcal S_{ji}^\prime := \partial_t (\mathcal S_{ji})$ and estimate $\left \| \mathcal S_{ji}^\prime \right\|_{\mathcal L \left(L^2((\partial D_j)_T),\,H^{1,1/2}((\partial D_i)_T)\right)}$ for $i\neq j$. Note that
\begin{eqnarray}\label{S_ij3}
&&\left\|\mathcal S_{ji}^\prime\right\|_{\mathcal L \left(L^2((\partial D_j)_T),\,H^{1,1/2}((\partial D_i)_T)\right)} \nonumber\\
&:=&\sup_{0\neq\varphi\in L^2((\partial D_j)_T)}\frac{\left\|\mathcal S_{ji}^\prime[\varphi]\right\|_{H^{1,1/2}((\partial D_i)_T)}}{\|\varphi\|_{L^2((\partial D_j)_T)}}\nonumber\\
&\leq& \sup_{0\neq\varphi\in L^2((\partial D_j)_T)}\frac{\left\|\mathcal S_{ji}^\prime[\varphi]\right\|_{L^2((\partial D_i)_T)} + \left\|\nabla_{\mathrm{tan}}\mathcal S_{ji}^\prime[\varphi]\right\|_{L^2((\partial D_i)_T)} + \left\|D_t^{1/2}\mathcal S_{ji}^\prime[\varphi]\right\|_{L^2((\partial D_i)_T)} }{\|\varphi\|_{L^2((\partial D_j)_T)}}.
\end{eqnarray}
We deduce that
\begin{eqnarray}\label{S_ij4}
\left\|\mathcal S_{ji}^\prime[\varphi]\right\|^2_{L^2((\partial D_i)_T)} &=& \int_0^T \int_{\partial D_i} \left| \int_0^t\int_{\partial D_j} \partial_t \Phi(x,\,t;\,y,\,\tau)\,\varphi(y,\,\tau)\,ds(y)d\tau \right|^2\,ds(x)dt\nonumber\\
&\lesssim& \|\varphi\|^2_{L^2((\partial D_j)_T)}\int_{\partial D_i}\int_{\partial D_j}\left(\int_0^T\int_0^t |\partial_t\Phi(x,\,t;\,y,\,\tau)|^2\, d\tau dt\right)ds(x)ds(y)\nonumber\\
&\lesssim& d_{ij}^{-8}\,\|\varphi\|^2_{L^2((\partial D_j)_T)}\,|\partial D_j|\,|\partial D_i|,
\end{eqnarray}
which implies that
\begin{equation}\label{S_ij5}
\left\|\mathcal S_{ji}^\prime[\varphi]\right\|_{L^2((\partial D_i)_T)} \lesssim \varepsilon^2\, d_{ij}^{-4}\,\|\varphi\|_{L^2((\partial D_j)_T)}.
\end{equation}
Analogously, we can also prove that
\begin{equation}\label{S_ij6}
\left\|\nabla_{\mathrm{tan}}\mathcal S_{ji}^\prime[\varphi]\right\|_{L^2((\partial D_i)_T)} \lesssim \varepsilon^2\, d_{ij}^{-5}\, \|\varphi\|_{L^2((\partial D_j)_T)}
\end{equation}
and
\begin{equation}\label{S_ij7}
\left\|D_t^{1/2}\mathcal S_{ji}^\prime[\varphi]\right\|_{L^2((\partial D_i)_T)} \lesssim \varepsilon^2\, d_{ij}^{-5}\, \|\varphi\|_{L^2((\partial D_j)_T)}.
\end{equation}
So we obtain from \eqref{S_ij3} that
\begin{equation}\label{S_ij2}
\left \| \mathcal S_{ji}^\prime \right\|_{\mathcal L \left(L^2((\partial D_j)_T),\,H^{1,1/2}((\partial D_i)_T)\right)} \lesssim \varepsilon^2\,d_{ij}^{-5},
\end{equation}
which leads to the following estimate:
\begin{equation}\label{S_ij8}
\Big\|\sum_{\substack{j=1\\j\neq i}}^M \partial_t \left(\mathcal S_{ji}[\sigma_j]\right)\Big\|_{H^{1,1/2}((\partial D_i)_T)} \lesssim \varepsilon^2 \sum_{\substack{j=1\\j\neq i}}^M d_{ij}^{-5} \|\sigma_j\|_{L^2((\partial D_j)_T)} \lesssim \varepsilon^2\,d^{-5} \sup_{j}\|\sigma_j\|_{L^2((\partial D_j)_T)}.
\end{equation}
As a consequence, we have proved that
\begin{equation}\label{est_sigma_t0}
\|\partial_t\sigma_i\|_{L^2((\partial D_i)_T)} \lesssim O(1) + \varepsilon\,d^{-5} \sup_{j} \|\sigma_j\|_{L^2((\partial D_j)_T)}.
\end{equation}
By repeatedly using this argument, we have
\begin{equation}\label{est_sigma_t}
\|\partial_t^n \sigma_i\|_{L^2((\partial D_i)_T)} \lesssim O(1) + \varepsilon\,d^{-3-2n} \sup_{j} \|\sigma_j\|_{L^2((\partial D_j)_T)}.
\end{equation}
Inserting \eqref{est_sigma_t} into \eqref{est_sigma_0}, we get
\begin{eqnarray}\label{est_sigma_final}
\| \sigma_i\|_{L^2((\partial D_i)_T)} &\lesssim& \sum_{\substack{j=1\\j\neq i}}^M \varepsilon d_{ji}^{-2} \|\sigma_j\|_{L^2((\partial D_j)_T)} + O(1) + \sum_{\substack{j=1\\j\neq i}}^M  O\left(\varepsilon^2\, d_{ji}^{-4}\,\|\sigma_j \|_{L^2((\partial D_j)_T)} \right) \nonumber\\
&& + O(\varepsilon\| \sigma_i \|_{L^2(\partial D_i)})+ O(\varepsilon) + O\left(\varepsilon^2\, d^{-5}\,\sup_{j} \|\sigma_j\|_{L^2((\partial D_j)_T)}\right).
\end{eqnarray}
Hence, we conclude that if $\varepsilon \sum_j\; d^{-2}_{ij} <1$, or equivalently, $1-2\beta-s/3\geq 0$, then
\begin{equation}\label{est_sigma_f}
\|\sigma_i\|_{L^2((\partial D_i)_T)} = O(1),\quad \|\partial_t \sigma_i\|_{L^2((\partial D_i)_T)} = O(1) \quad \mbox{ as } \varepsilon \ll 1.
\end{equation}

We state this results in the following proposition:
\begin{prop}\label{a-priori-estimates}
Under the following condition on the distribution of the small cavities
\begin{equation}\label{condi-a-priori-estimates}
\varepsilon\; \max_{1\leq i\leq M}\sum_{j\neq i}\; d^{-2}_{ij} <1
\end{equation}
the solution of the system of integral equation \eqref{ie_bvp_M_1} has the following estimates
\begin{equation}\label{the-apr-estimates}
\|\sigma_i\|_{L^2((\partial D_i)_T)} = O(1),\quad \|\partial_t \sigma_i\|_{L^2((\partial D_i)_T)} = O(1) \quad \mbox{ as } \varepsilon \ll 1.
\end{equation}
\end{prop}

\subsection{Invertibility of the algebraic system}

For $j=1,\,2,\,\cdots,\,M$, we define
\begin{equation}\label{def_Ci}
C_j := \int_{\partial D_j} S_{\partial D_j}^{-1}[1](y)\,ds(y).
\end{equation}
We state the following system of integral equations
\begin{equation}\label{MS7}
q_i(t) + \sum_{\substack{j=1\\j\neq i}}^M \int_0^t C_j\, \Phi(z_i,\,t;\,z_j,\,\tau)\,q_j(\tau)\,d\tau = f_i(t), \quad i=1,\,2,\,\cdots,\, M.
\end{equation}
This system is naturally linked to the proof of our main results in Section \ref{End-proof-Main-asymptotic}. Here, we show the invertibility of the system \eqref{MS7} and estimate $\sum_{i=1}^M \|q_i\|^2_{L^2(0,\,T)}$.

\begin{theorem}\label{sys_sol}
If
\begin{equation}\label{MS8}
C \max_{1\leq i \leq M}\sum_{j\neq i}\vert z_i-z_j \vert^{-2} <1.
\end{equation}
with $C:=\max_{1\leq j \leq M}C_j$, then the system \eqref{MS7} is uniquely solvable. Moreover, we have the estimate
\begin{eqnarray}\label{MS_q}
\left(\sum_{i=1}^M \|q_i\|_{L^2}^2\right)^{1/2} \leq \left(1- C \max_{1\leq i \leq M}\sum_{j\neq i}\vert z_i-z_j \vert^{-2}\right)^{-1}  \Bigg\{\sum_{i=1}^M \Vert f_i\Vert^2_{L^2(0, T)} \Bigg\}^{1/2}.
\end{eqnarray}
\end{theorem}

{\bf Proof.} Observe from \eqref{MS7} that
\begin{eqnarray}\label{MS9}
\sum_{i=1}^M \int_0^T q_i^2(t)\,dt + \sum_{i=1}^M \sum_{\substack{j=1\\j\neq i}}^M C_j \int_0^T \int_0^t \Phi(z_i,\,t;\,z_j,\,\tau) \,q_j(\tau)\,q_i(t)\,d\tau dt =  \sum_{i=1}^M \int_0^T f_i(t)\,q_i(t)\,dt.
\end{eqnarray}

Note that
\begin{eqnarray}\label{MS11}
\int_0^T \int_0^t \left|\Phi(z_i,\,t;\,z_j,\,\tau)\, q_j(\tau)\,q_i(t) \right| \,d\tau dt
&\leq & \int_0^T \left( \int_0^t |\Phi(z_i,\,t;\,z_j,\,\tau)|^2\,d\tau \right)^{1/2} \|q_j\|_{L^2} \,|q_i(t)|\,dt \nonumber\\
&\leq & \left(\int_0^T  \int_0^t |\Phi(z_i,\,t;\,z_j,\,\tau)|^2\,d\tau dt \right)^{1/2}\,\|q_j\|_{L^2}\,
\|q_i\|_{L^2}.
\end{eqnarray}
Using Lemma \ref{int_phi}, we deduce from \eqref{MS9} that
\begin{eqnarray}\label{MS12}
\sum_{i=1}^M \|q_i\|_{L^2}^2 - C  \sum_{i=1}^M \sum_{\substack{j=1\\j\neq i}}^M |z_i-z_j|^{-2}\,\|q_i\|_{L^2}^2   \leq \left( \sum_{i=1}^M \|f_i\|_{L^2}^2 \right)^{1/2}\, \left( \sum_{i=1}^M \|q_i\|_{L^2}^2 \right)^{1/2}.
\end{eqnarray}
Then we have
\begin{equation*}
\left(1- C\max_{i}\sum_{\substack{j=1\\j\neq i}}^M |z_i-z_j|^{-2}\right)\,\sum_{i=1}^M \|q_i\|_{L^2}^2 \leq \left( \sum_{i=1}^M \|f_i\|_{L^2}^2 \right)^{1/2}\, \left( \sum_{i=1}^M \|q_i\|_{L^2}^2 \right)^{1/2},
\end{equation*}
and therefore
\begin{equation}\label{MS14}
\left(1- C \max_{i}\sum_{\substack{j=1\\j\neq i}}^M |z_i-z_j|^{-2}\right) \, \left(\sum_{i=1}^M \|q_i\|_{L^2}^2\right)^{1/2} \leq \left( \sum_{i=1}^M \|f_i\|_{L^2}^2 \right)^{1/2}.
\end{equation}
Thus, the unique solvability of \eqref{MS7} and the estimate \eqref{MS_q} follow from the condition \eqref{MS8}. The proof is now complete. \hfill $\Box$

\subsection{End of the proof of Theorem \ref{Main-asymptotic}}\label{End-proof-Main-asymptotic}

In the sequel, we derive the asymptotic formula for the solution to \eqref{bvp} in the case of multiple cavities. For $x\in\partial D_i$ with $i\neq j$, we derive that
\begin{eqnarray}\label{ie_bvp_M00}
&&\int_0^t\int_{\partial D_j} \Phi(x,\,t;\,y,\,\tau)\,\sigma_j(y,\,\tau)\,ds(y)d\tau  \nonumber\\
&=& \int_0^t \Phi(x,\,t;\,z_j,\,\tau)\,\left( \int_{\partial D_j} \sigma_j(y,\,\tau)\,ds(y) \right)d\tau \nonumber \\
&&+ \int_0^t\int_{\partial D_j} \left[ \Phi(x,\,t;\,y,\,\tau) - \Phi(x,\,t;\,z_j,\,\tau) \right]\,\sigma_j(y,\,\tau)\,ds(y)d\tau  \nonumber\\
&=& \int_0^t \Phi(x,\,t;\,z_j,\,\tau)\,\left( \int_{\partial D_j} \sigma_j(y,\,\tau)\,ds(y) \right)d\tau + O(d_{ij}^{-3}\,\varepsilon^2) \|\sigma_j\|_{L^2((\partial D_j)_T)} \nonumber\\
&=& \int_0^t \Phi(z_i,\,t;\,z_j,\,\tau)\,\left( \int_{\partial D_j} \sigma_j(y,\,\tau)\,ds(y) \right)d\tau  \nonumber\\
&& + \int_0^t \left[\Phi(x,\,t;\,z_j,\,\tau) - \Phi(z_i,\,t;\,z_j,\,\tau)\right]\,\left( \int_{\partial D_j} \sigma_j(y,\,\tau)\,ds(y) \right)d\tau + O(d_{ij}^{-3}\varepsilon^2)\,\|\sigma_j\|_{L^2((\partial D_j)_T)} \nonumber\\
&=& \int_0^t \Phi(z_i,\,t;\,z_j,\,\tau)\,\left( \int_{\partial D_j} \sigma_j(y,\,\tau)\,ds(y) \right)d\tau + O(d_{ij}^{-3}\,\varepsilon^2)\, \|\sigma_j\|_{L^2((\partial D_j)_T)},
\end{eqnarray}
since, by using Taylor's expansion and \eqref{int_dphi_e}, we have
\begin{eqnarray}\label{ie_bvp_M01}
&&\left|\int_0^t \int_{\partial D_j} \left[\Phi(x,\,t;\,y,\,\tau) - \Phi(x,\,t;\,z_j,\,\tau)\right]\, \sigma_j(y,\,\tau)\,ds(y)d\tau\right| \nonumber\\
&\leq & \left(\int_0^t\int_{\partial D_j} |\nabla_y \Phi(x,\,t;\,z_j^*,\,\tau)|^2\,|y-z_j|^2\,ds(y)d\tau\right)^{1/2}\,\left(\int_0^t \int_{\partial D_j} \sigma_j^2(y,\,\tau)\,ds(y)d\tau \right)^{1/2}\nonumber\\
&\lesssim & \varepsilon^2\, d_{ij}^{-3} \, \|\sigma_j\|_{L^2((\partial D_j)_T)}
\end{eqnarray}
and
\begin{eqnarray}\label{ie_bvp_M02}
&&\left|\int_0^t \left[\Phi(x,\,t;\,z_j,\,\tau) - \Phi(z_i,\,t;\,z_j,\,\tau)\right]\,\left( \int_{\partial D_j} \sigma_j(y,\,\tau)\,ds(y) \right)d\tau\right| \nonumber\\
&\leq & \left(\int_0^t |\nabla_x \Phi(z_i^*,\,t;\,z_j,\,\tau)|^2\,|x-z_i|^2\,d\tau\right)^{1/2}\,\left[\int_0^t \left( \int_{\partial D_j} \sigma_j(y,\,\tau)\,ds(y) \right)^2\,d\tau \right]^{1/2}\nonumber\\
&\lesssim & \varepsilon^2\, d_{ij}^{-3} \, \|\sigma_j\|_{L^2((\partial D_j)_T)}
\end{eqnarray}
where $z_j^*=z_j + \theta_1\, (y-z_j),\,y\in D_j,\,0<\theta_1<1$ and $z_i^*=z_i + \theta_2\, (x-z_i),\,x\in D_i,0<\theta_2<1$.
Hence, the boundary integral system \eqref{ie_bvp_M} can be rewritten as
\begin{eqnarray}\label{ie_bvp_M1}
&&\int_0^t \int_{\partial D_i} \Phi(x,\,t;\,y,\,\tau)\,\sigma_i(y,\,\tau)\,ds(y)d\tau + \sum_{\substack{j=1\\j\neq i}}^M \int_0^t  \Phi(z_i,\,t;\,z_j,\,\tau)\,\left(\int_{\partial D_j} \sigma_j(y,\,\tau)\,ds(y)\right) d\tau \nonumber\\
&=& \Phi(z_i,\,t;\,z^*,\,0) + O\left( \sum_{\substack{j=1\\j\neq i}}^M d_{ij}^{-3}\,\varepsilon^2\right) + O(\varepsilon), \quad x\in\partial D_i,\; i=1,\,2,\,\cdots,\,M.
\end{eqnarray}

Recall that
\begin{equation*}
\varphi_i(x,\,y,\,t) = \int_0^t \frac{|x-y|}{2\sqrt{\pi}(t-\tau)^{3/2}}\,\exp\left( -\frac{|x-y|^2}{4(t-\tau)} \right)\, \sigma_i(y,\,\tau)\,d\tau, \quad x\in\mathbb R^3,\;y\in\partial D_i,\;t\in (0,\,T].
\end{equation*}
Then, by Lemma \ref{I}, we have
\begin{eqnarray}\label{ie_bvp_M2}
&&\int_0^t \int_{\partial D_i} \Phi(x,\,t;\,y,\,\tau)\,\sigma_i(y,\,\tau)\,ds(y)d\tau\nonumber\\
&=& \int_{\partial D_i}\frac{1}{4\pi |x-y|}\, \sigma_i(y,\,t)\,ds(y) + \int_{\partial D_i}\frac{1}{4\pi |x-y|} \left[ \varphi_i(x,\,y,\,t) - \sigma_i(y,\,t)\right]\,ds(y) \nonumber\\
&=& \int_{\partial D_i}\frac{1}{4\pi |x-y|}\, \sigma_i(y,\,t)\,ds(y) +  O\left(\int_{\partial D} \|\partial_t \sigma_i \|_{L^2(0,\,t)}\right) \nonumber\\
&=& \int_{\partial D_i}\frac{1}{4\pi |x-y|}\, \sigma_i(y,\,t)\,ds(y) +  O\left(\varepsilon\right), \quad x\in \partial D_i,\;t\in(0,\,T].
\end{eqnarray}
It follows from \eqref{ie_bvp_M1} and \eqref{ie_bvp_M2} that
\begin{eqnarray}\label{ie_bvp_M3}
&&\int_{\partial D_i}\frac{1}{4\pi |x-y|} \sigma_i(y,\,t)\,ds(y) + \sum_{\substack{j=1\\j\neq i}}^M \int_0^t  \Phi(z_i,\,t;\,z_j,\,\tau)\,\left(\int_{\partial D_j} \sigma_j(y,\,\tau)\,ds(y)\right) d\tau  \nonumber\\
&=& \Phi(z_i,\,t;\,z^*,\,0) + O\left( \sum_{\substack{j=1\\j\neq i}}^M d_{ij}^{-3}\,\varepsilon^2 \right) + O(\varepsilon).
\end{eqnarray}
Hence, we deduce that
\begin{eqnarray}\label{ie_bvp_M4}
&&\int_{\partial D_i} \sigma_i(y,\,t)\,ds(y) + \left( \int_{\partial D_i} S_{\partial D_i}^{-1}[1](y)\,ds(y)\right) \sum_{\substack{j=1\\j\neq i}}^M \int_0^t \Phi(z_i,\,t;\,z_j,\,\tau) \left(\int_{\partial D_j} \sigma_j(y,\,\tau)\,ds(y)\right)d\tau\nonumber\\
&=&\Phi(z_i,\,t;\,z^*,\,0) \int_{\partial D_i} S_{\partial D_i}^{-1}[1](y)\,ds(y) + E_i\,\int_{\partial D_i} S_{\partial D_i}^{-1}[1](y)\,ds(y),
\end{eqnarray}
where
\begin{equation*}
E_i := O\left( \sum_{\substack{j=1\\j\neq i}}^M d_{ij}^{-3}\,\varepsilon^2 \right) + O(\varepsilon)
=O\left(\frac{\varepsilon^2\,|\ln \varepsilon|}{d^3}\right) + O(\varepsilon)
=O\left(\varepsilon^{2-3\beta}\,|\ln \varepsilon|\right) + O(\varepsilon).
\end{equation*}
Then we have the following system:
\begin{equation}\label{ie_bvp_M5}
\frac{q_i(t)}{C_i} + \sum_{\substack{j=1\\j\neq i}}^M \int_0^t C_j\, \Phi(z_i,\,t;\,z_j,\,\tau)\,\left(\frac{q_j(\tau)}{C_j}\right) d\tau = \Phi(z_i,\,t;\,z^*,\,0) + O\left(\varepsilon^{2-3\beta}\,|\ln \varepsilon|\right) + O(\varepsilon),\quad i=1,\,2,\,\cdots,\,M.
\end{equation}

We are now in a position to state our main result of this section.
\begin{theorem}\label{th_asym2}
For $x\in \mathbb R^3\setminus\overline D$ and $t\in(0,\,T]$, the solution to \eqref{bvp} has the following asymptotic expansion:
\begin{equation}\label{asym_M}
u(x,\,t) = \sum_{i=1}^M C_i \int_0^t \Phi(x,\,t;\,z_i,\,\tau)\,\alpha_i(\tau)\, d\tau
+ O(\varepsilon^{3-3\beta-s} |\ln \varepsilon|) + O(\varepsilon^{2-s}) \quad \textrm{as}\; \varepsilon \to 0 \end{equation}
under the following condition on the distribution of the cavities
\begin{equation}\label{condi-a-priori-estimates-2}
\varepsilon\; \max_{1\leq i \leq M}\sum_{j\neq i}\; d^{-2}_{ij} <1,
\end{equation}
or $1-2\beta -\frac{s}{3}\geq 0$,  where the constants $C_i$'s are defined by \eqref{def_Ci} and $\left\{ \alpha_i(t) \right\}_{i=1}^M$ is the unique solution of the linear system
\begin{equation}\label{ie_bvp_M13}
\alpha_i(t) + \sum_{\substack{j=1\\j\neq i}}^M \int_0^t C_j \,\Phi(z_i,\,t;\,z_j,\,\tau)\,\alpha_j(\tau)\, d\tau = \Phi(z_i,\,t;\,z^*,\,0),\quad i=1,\,2,\,\cdots,\,M.
\end{equation}
\end{theorem}

{\bf Proof.} Define
\begin{equation*}
w_i(t) := q_i(t)/C_i - \alpha_i(t),\quad i=1,\,2,\,\cdots,\,M.
\end{equation*}
Then, by \eqref{ie_bvp_M5} and \eqref{ie_bvp_M13}, we have
\begin{equation}\label{ie_bvp_M14}
w_i(t) + \sum_{\substack{j=1\\j\neq i}}^M \int_0^t C_j\, \Phi(z_i,\,t;\,z_j,\,\tau)\,w_j(\tau)\, d\tau = O\left(\varepsilon^{2-3\beta}\,|\ln \varepsilon|\right) + O(\varepsilon):=\overline{E}_i,\quad i=1,\,2,\,\cdots,\,M.
\end{equation}
From Theorem \ref{sys_sol}, we obtain that
\begin{equation}\label{ie_bvp_M15}
\sum_{i=1}^M \|w_i\|_{L^2}^2 \lesssim\sum_{i=1}^M \|\overline{E}_i \|_{L^2}^2 =: E_1 \sim M (\varepsilon^{2-3\beta}\,
|\ln \varepsilon|+\varepsilon)^2,
\end{equation}
Thus, we derive for $x\in \mathbb R^3\setminus\overline D$ and $t\in (0,\,T]$ that
\begin{eqnarray}\label{ie_bvp_M17}
u(x,\,t) &=& \sum_{i=1}^M\int_0^t \int_{\partial D_i} \Phi(x,\,t;\,y,\,\tau)\,\sigma_i(y,\,\tau)\,ds(y)d\tau \nonumber\\
&=& \sum_{i=1}^M\int_0^t \int_{\partial D_i} \Phi(x,\,t;\,z_i,\,\tau)\,\sigma_i(y,\,\tau)\,ds(y)d\tau \nonumber\\
&& + \sum_{i=1}^M\int_0^t \int_{\partial D_i} \left[\Phi(x,\,t;\,y,\,\tau) - \Phi(x,\,t;\,z_i,\,\tau) \right]\,\sigma_i(y,\,\tau)\,ds(y)d\tau \nonumber\\
&=& \sum_{i=1}^M \int_0^t \Phi(x,\,t;\,z_i,\,\tau)\,q_i(\tau)\, d\tau + O(\varepsilon^2)\sum_{i=1}^M\Vert \sigma_i\Vert_{L^2(\partial D_i)_T} \nonumber\\
&=& \sum_{i=1}^M C_i \int_0^t \Phi(x,\,t;\,z_i,\,\tau)\,\alpha_i(\tau)\, d\tau + O(C(ME_1)^{\frac{1}{2}}) + O(M\varepsilon^2) \nonumber\\
&=& \sum_{i=1}^M C_i \int_0^t \Phi(x,\,t;\,z_i,\,\tau)\,\alpha_i(\tau)\, d\tau + O(\varepsilon^{3-3\beta-s} |\ln \varepsilon|) + O(\varepsilon^{2-s})  \quad \textrm{for}\; M\sim \varepsilon^{-s},\; d \sim \varepsilon^{\beta}.
\end{eqnarray}
This completes the proof. \hfill $\Box$

\section{Proof of Theorem \ref{Main-effective}}\label{effective-medium}
\setcounter{equation}{0}

Let $\Omega$ be a bounded domain containing the cavities $D_j,\,j=1,\,2,\,\cdots,\,M$. We divide $\Omega$ into $[ a^{-1}]$ periodically subdomains $\Omega_j,\,j=1,\,2,\,\cdots,\,[a^{-1}]$ such that $\Omega_j$'s are disjoint and each $\Omega_j$ contains one single cavity $D_j$ and has a volume $a$. We also assume that the cavities $D_j,\,j=1,\,2,\,\cdots,\,M$ have the same shape. This means that $C_i=C_j$ for $i,\,j=1,\,2,\,\cdots,\,M$. Define $$C:=C_j=\overline C\; a,$$ where $\overline C$ is the scaled value of $C_j$.

As $\Omega$ can have an arbitrary shape, the set of the cubes intersecting $\partial \Omega$ is not empty (unless if $\Omega$ has a simple shape as a cube). Later in our analysis, we will need the estimate of the volume of this set. Since each $\Omega_j$ has volume of the order $a$, and then its maximum radius is of the order $a^{\frac{1}{3}}$, then the intersecting surfaces with $\partial \Omega$ has an area of the order $a^{\frac{2}{3}}$. As the area of $\partial \Omega$ is of the order one, we conclude that the number of such cubes will not exceed the order $a^{-\frac{2}{3}}$. Hence the volume of this set will not exceed the order $a^{-\frac{2}{3}}a=a^{\frac{1}{3}}$, as $a \rightarrow 0$.

We consider the integral equation
\begin{equation}\label{med_1}
v(x,\,t) + \int_0^t \int_\Omega \overline C\,\Phi(x,\,t;\,z,\,\tau)\,v(z,\,\tau)\,dzd\tau = \Phi(x,\,t;\,z^*,\,0),\quad (x,\,t)\in\Omega_T,\; z^*\not\in\overline\Omega.
\end{equation}
The unique solvability can be proved as follows.

\begin{lemma}\label{med_sol}
The integral equation \eqref{med_1} is uniquely solvable in $L^2(\Omega_T)$.
\end{lemma}

{\bf Proof.} Due to the estimate \eqref{fun_est1}, the volume potential operator $\mathcal V$ defined by
\begin{equation*}
\mathcal V[\varphi](x,\,t):= \int_0^t \int_\Omega \Phi(x,\,t;\,z,\,\tau)\,\varphi(z,\,\tau)\,dzd\tau
\end{equation*}
has a weakly singular kernel and is bounded from $L^2(\Omega_T)$ into $L^2(\Omega_T)$. We note that the operator norm of $\mathcal V$ goes to zero as $T\to 0$. Thus the inverse of $I+\mathcal V$ can be written as the series $\sum_{j=0}^\infty (-\mathcal V)^j$ for small $T$. The result for large $T$ follows by an iteration argument. \hfill $\Box$

\bigskip
Define
\begin{equation}\label{def_V}
V(x,\,t):=
\begin{cases}
v(x,\,t) & \mathrm{in}\;\Omega_T,\\
\Phi(x,\,t;\,z^*,\,0) - \displaystyle\int_0^t \int_\Omega \overline C\,\Phi(x,\,t;\,z,\,\tau)\,v(z,\,\tau)\,dzd\tau & \mathrm{in}\; (\mathbb R^3\setminus\overline \Omega)_T.
\end{cases}
\end{equation}
Then we have the following result:
\begin{lemma}\label{V_v}
$v$ is the solution of \eqref{med_1} if and only if $V$ is a solution of $(\partial_t - \Delta + \overline C \chi_\Omega) V = \delta(x-z^*)\,\delta(t)$ in $\mathbb R^3$.
\end{lemma}

Set
\begin{equation*}
W(x,\,t) := \Phi(x,\,t;\,z^*,\,0) - V(x,\,t).
\end{equation*}
Then we obtain that
\begin{equation}\label{W-1}
\begin{cases}
(\partial_t - \Delta + \overline C \chi_\Omega)W = \overline C \chi_\Omega \Phi(x,\,t;\,z^*,\,0) & \mathrm{in}\;\mathbb R^3\times (0,\,T), \\
W(x,\,0)=0 & \mathrm{in}\; \mathbb R^3, \\
|W(x,\,t)|\leq C_0 \exp(b|x|^2) & \mathrm{as}\; |x|\to +\infty.
\end{cases}
\end{equation}
As $z^*$ is outside $\Omega$, the right hand side of the first equation in (\ref{W-1}) is smooth.
Then the solution $W$ is in $ C\left([0,\,T];\,H_{loc}^1(\mathbb R^3) \right)$; see \cite{Is} for instance. By Sobolev embedding, we deduce that $W\in C\left([0,\,T];\,L^p(\Omega) \right)$ for any $p<3$.

From \eqref{med_1}, we see
$$
\vert v(x,\,t)\vert \lesssim \int_0^t\int_\Omega |\Phi(x,\,t;\,z,\,\tau)| |v(z, \tau) | \,dzd\tau+|\Phi(x,\,t;\,z^*,\,0)|
$$
and hence
$$
\vert v(x,\,t)\vert \lesssim \left(\int_0^t\int_\Omega |\Phi(x,\,t;\,z,\,\tau)|^q\,dzd\tau\right)^\frac{1}{q} \Vert v \Vert_{C\left([0,\,T];\,L^p(\Omega) \right)}+O(1), \quad \mbox{ with } \frac{1}{p}+\frac{1}{q}=1.
$$
By the singularity estimate \eqref{fun_est1}, we have
$$
|\Phi(x,\,t;\,z,\,\tau)|^q \lesssim \frac{1}{(t-\tau)^{\mu q}} \frac{1}{\vert x-z\vert^{(3-2\mu)q}}
$$
and this function is integrable in $\Omega \times (0, \,T)$ if $\mu q<1$ and $(3-2\mu)q<3$.
As $p<3$, then $q>\frac{3}{2}$. Choosing $\mu$ smaller but near $\frac{2}{3}$, then these two conditions on $q$ are satisfied.
Hence $v$ is in $C\left([0,\,T]; \,L^\infty(\Omega) \right)$. In addition, from \eqref{med_1}, we get
\begin{equation}\label{med_3}
\partial_{x_j} v(x,\,t) = - \int_0^t\int_\Omega \overline C\,\partial_{x_j}\Phi(x,\,t;\,z,\,\tau)\,v(z,\,\tau)\,dzd\tau + \partial_{x_j}\Phi(x,\,t;\,z^*,\,0), \quad (x,\,t)\in\Omega_T,\; z^*\not\in\Omega,
\end{equation}
and then
\begin{equation}\label{med_4}
|\partial_{x_j}v(x,\,t) |
\lesssim \left(\int_0^t\int_\Omega |\partial_{x_j}\Phi(x,\,t;\,z,\,\tau)|\,dzd\tau\right)\, \|v \|_{L^\infty\left((0,\,T);\,L^\infty(\Omega) \right)} + |\partial_{x_j}\Phi(x,\,t;\,z^*,\,0)|.
\end{equation}
By the singularity estimate \eqref{fun_est2} with $\mu>1/2$, we see
\begin{equation*}
\int_0^t\int_\Omega |\partial_{x_j}\Phi(x,\,t;\,z,\,\tau)|\,dz = O(1).
\end{equation*}
This means that $\partial_{x_j} v\in C\left([0,\,T];\,L^\infty(\Omega) \right)$. So we obtain that $v\in C\left([0,\,T];\, W^{1,\infty}(\Omega)\right)$.

\bigskip
We now rewrite the integral equation \eqref{med_1} at $x=z_l$ for $1\leq l \leq M$ as
\begin{equation}\label{med_1_1}
v(z_l,\,t) + \int _0^t \sum_{\substack{j=1\\j\neq l}}^{[a^{-1}]} \int_{\Omega_j} \overline C\,\Phi(z_l,\,t;\,z,\,\tau)\,v(z,\,\tau)\,dzd\tau = \Phi(z_l,\,t;\,z^*,\,0) + \mathcal A +\mathcal A_l,
\end{equation}
or
\begin{equation}\label{med_1_2}
v(z_l,\,t) + \sum_{\substack{j=1\\j\neq l}}^{M} \overline C\, {a} \int_0^t \Phi(z_l,\,t;\,z_j,\,\tau)\,v(z_j,\,\tau)\,dzd\tau = \Phi(z_l,\,t;\,z^*,\,0) + \mathcal A + \mathcal A_l + \mathcal B_l
\end{equation}
with
\begin{eqnarray*}
\mathcal A &:= & - \int_0^t \int_{\Omega \setminus(\cup^{[a^{-1}]}_{j=1}\Omega_j)} \overline C\,\Phi(z_l,\,t;\,z,\,\tau)\,v(z,\,\tau)\,dzd\tau,\\
\mathcal A_l & := & - \int_0^t\int_{\Omega_l} \overline C\,\Phi(z_l,\,t;\,z,\,\tau)\,v(z,\,\tau)\, dzd\tau,\\
\mathcal B_l & := & - \sum_{\substack{j=1\\j\neq l}}^{[a^{-1}]}\overline C \int_0^t \int_{\Omega_j} \,\Phi(z_l,\,t;\,z,\,\tau)\,v(z,\,\tau)\,dz d \tau +  \sum_{\substack{j=1\\j\neq l}}^{M}\, \overline C\, {a}\, \int^t_0\Phi(z_l,\,t;\,z_j,\,\tau)\,v(z_j,\,\tau) \,d\tau.
\end{eqnarray*}

As $v\in L^{\infty}(\Omega_T)$, then by (\ref{fun_est1}), we have
$$\mathcal A_l \sim \int^t_0 \int_{\Omega_l} \vert \Phi(z_l,\,t;\,z,\,\tau)\vert \,dz d\tau = O\left( \int_{\Omega_l} \vert z-z_l\vert^{2\mu-3}\, dz\right)$$ for $0<\mu<1$, and hence, by a scaling, we derive the estimate
\begin{equation}\label{est_A_l}
\mathcal A_l =O\left(\varepsilon^{\frac{2\mu}{3}}\right)\quad \mathrm{as}\; \varepsilon \ll 1.
\end{equation}

\medskip
Let us estimate $\mathcal B_l$. As $\vert \Omega_l\vert= a $, we have
$$
\mathcal B_l=- \sum_{\substack{j=1\\j\neq l}}^{[a^{-1}]} \overline C \int_0^t \int_{\Omega_j} \left[\Phi(z_l,\,t;\,z,\,\tau)\,v(z,\,\tau) - \Phi(z_l,\,t;\,z_j,\,\tau)\,v(z_j,\,\tau)\right] \,dz d \tau.
$$
We write the above integrand as
\begin{eqnarray*}
&&\Phi(z_l,\,t;\,z,\,\tau)\,v(z,\,\tau) - \Phi(z_l,\,t;\,z_j,\,\tau)\,v(z_j,\,\tau)\\ & = &
[\Phi(z_l,\,t;\,z,\,\tau)-\Phi(z_l,\,t;\,z_j,\,\tau)]v(z,\,\tau) +  \Phi(z_l,\,t;\,z_j,\,\tau)[v(z,\,\tau)-v(z_j,\,\tau)].
\end{eqnarray*}
Then we see
$$
\mathcal B_l = O\Big(\sum_{\substack{j=1\\j\neq l}}^{[a^{-1}]} \int^t_0 \int_{\Omega_j} \left[\vert \nabla_z\Phi(z_l,\,t;\,z_j,\,\tau)\vert \, \vert z-z_j\vert +\vert \Phi(z_l,\,t;\,z_j,\,\tau)\vert\, \vert z-z_j\vert\right]\,dz d\tau\Big).
$$
But
\begin{eqnarray*}
\int_{\Omega_j} \vert\Phi(z_l,\,t;\,z_j,\,\tau)\vert\,\vert z-z_j\vert\, dz &=& O\left((t-\tau)^{-\mu}\; \vert z_l-z_j\vert^{-3+2\mu}\right)
\int_{\Omega_j}\vert z-z_j\vert\, dz\\
&=&O\left((t-\tau)^{-\mu}\; \vert z_l-z_j\vert^{-3+2\mu}\right) a^{\frac{4}{3}}
\end{eqnarray*}
with $0<\mu<1$, and similarly
$$
\int_{\Omega_j} \vert \nabla_z \Phi(z_l,\,t;\,z_j,\,\tau)\vert\,\vert z-z_j\vert\, dz =O\left((t-\tau)^{-\gamma}\; \vert z_l-z_j\vert^{-4+2\gamma}\right)\, a^{\frac{4}{3}}
$$
with $0< \gamma<1$. Hence, choosing $\mu = \gamma \in [1/2,\,1)$, we get
\begin{equation}\label{est_B_l}
\mathcal B_l=O\Big(\sum_{\substack{j=1\\j\neq l}}^{[a^{-1}]} \vert z_l-z_j\vert^{-4+2\gamma}\Big)\, a^{\frac{4}{3}}=O\left(d^{-3}\, a^{\frac{4}{3}}\right)=O\left(a^{\frac{1}{3}}\right).
\end{equation}

To estimate the term $\mathcal A$, following \cite{A-C-C-S-18, C-M-S-17}, we distinguish between the following two cases:
\begin{itemize}
\item[(a)] The point $z_m$ is away from the boundary $\partial \Omega  $ and so $\Phi(z_m,\, t;\, z,\,\tau) $ is bounded in $z$ near the boundary.

\item[(b)] The point $z_m$ is located near one of the $\Omega_j $'s touching the boundary $\partial \Omega $. In this case, we split the estimate into two parts. By $ N_m$ we denote the part that involves $\Omega_j $'s close to $z_m $, and we denote the remaining part by $F_m $. The integral over $F_m $ can be estimated in a manner similar to the case $(a)$ discussed above. Also note that $F_m \subset \Omega\setminus \cup_{j=1}^{[a^{-1}]} \Omega_j $ and so $\mathrm{Vol}\,(F_m) $ is of the order $a^{\frac{1}{3}} $ as $a \rightarrow 0 $.

\medskip
To estimate the integral over $N_m $, we observe that owing to the fact $a$ is small, the $\Omega_j $'s  close to $z_m $ are located near a small region of the boundary $\partial \Omega $. Since we assume that the boundary is smooth enough, this region can be assumed to be flat. We now divide this layer into concentric layers as in the estimate of $\mathcal B_l$. In this case, we have at most $(2n+1)^{2} $ cubes intersecting the surface, for $n=0,\,\dots,\,[a^{-\frac{1}{3}}] $. So the number of cavities in the $n^{th} $ layer $(n \neq 0)$ will be at most $[(2n+1)^{2}-(2n-1)^{2}] $ and their distance from $\Omega_m $ is at least $n\left(a^{\frac{1}{3}}-\frac{a}{2} \right)$.
\end{itemize}
Similar as in Lemma \ref{int_dphi}, we have $\int^t_0 \Phi(x,\, t;\, z,\, \tau)\, d\tau =O(\vert x-z\vert^{-1})$. Therefore we can write
\begin{equation*}
\begin{aligned}
\vert \mathcal A \vert&=\left\vert \int^t_0 \int_{\Omega \setminus (\cup_{j=1 }^{[a^{-1}]} \Omega_j)} \overline C\,  \Phi(z_m,\, t;\, z,\, \tau)  \, v(z,\, \tau) \, dz d\tau \right\vert\\
 &=\left\vert \int^t_0 \int_{N_m} \overline C\, \Phi(z_m,\, t;\, z,\, \tau)\,  v(z,\,\tau)\, dz d\tau \right\vert + \left\vert \int^t_0 \int_{F_m} \overline C\, \Phi(z_m,\, t;\, z,\, \tau)\,  v(z,\,\tau) \, dzd\tau \right\vert \\
&\leq \sum_{l=1}^{[a^{-\frac{2}{3}}]}\overline C\,  \Vert v\Vert_{L^{\infty}(\Omega_T)} \,\mathrm{Vol}\,(\Omega_l)\, \frac{1}{d_{ml}} +  \overline C\, \Vert{\Phi(z_m,\,t;\,\cdot,\,\tau)}\,\Vert_{L^{\infty}((F_m)_T)} \, \Vert v \Vert_{L^{\infty}(\Omega_T)}\, \mathrm{Vol}\,(F_m) \\
&\leq O\Big(a \sum_{l=1}^{[a^{-\frac{1}{3}}]} \frac{1}{d_{ml}} + \tilde C\, a^{\frac{1}{3}}\Big) \\
&\leq O\Big( a \sum_{l=1}^{[a^{-\frac{1}{3}}]} \Big[(2n+1)^2-(2n-1)^2 \Big]\, \frac{1}{n \Big( a^{\frac{1}{3}}-\frac{a}{2} \Big)} + \tilde C\,a^{\frac{1}{3}}\Big)\\
&= O\Big( a \,O(a^{-\frac{2}{3}}) + O(a^{\frac{1}{3}})\Big),
\end{aligned}
\label{est-vol-2}
\end{equation*}
and hence
\begin{equation}\label{est-vol-2a}
\vert \mathcal A \vert = O\left(a^{\frac{1}{3}}\right).
\end{equation}

Gathering the estimates \eqref{est_A_l}, \eqref{est_B_l} and \eqref{est-vol-2a}, taking $\mu = 1/2$, we have
$$
\sum_l \left(\vert \mathcal A \vert^2 +\vert \mathcal A_l \vert^2 +\vert \mathcal B_l\vert^2\right)=
O\left(Ma^{\frac{2}{3}}+ Ma^{\frac{4\mu}{3}}\right)=O\left(a^{-\frac{1}{3}}\right).
$$
Using the invertibility property and the estimate \eqref{MS_q} for the algebraic system \eqref{MS7}, we deduce the following estimate:
\begin{equation}\label{med_8}
\sum_{j=1}^M \|\alpha_j(t) - v(z_j,\,t) \|_{L^2(0,\,T)}^2 = O(a^{-\frac{1}{3}}) \quad \mathrm{as}\; a \to 0.
\end{equation}

Let $x$ be away from $\Omega \cup \{z^*\}$. Then we recall that
\begin{equation*}
W(x,\,t) := \Phi(x,\,t;\,z^*,\,0) - V(x,\,t) = \int_0^t \int_\Omega \overline C\, \Phi(x,\,t;\,z,\,\tau)\,v(z,\,\tau)\,dzd\tau,
\end{equation*}
and rewrite it as
\begin{equation*}
W(x,\,t) = \sum_{j=1}^{[a^{-1}]} |\Omega_j| \int_0^t \overline C \,\Phi(x,\,t;\,z_j,\,\tau)\,v(z_j,\,\tau)\,d\tau + \mathcal C
\end{equation*}
with
\begin{equation*}
\mathcal C = \sum_{j=1}^{[a^{-1}]} \int_0^t \int_{\Omega_j}\overline C\,\left[\Phi(x,\,t;\,z,\,\tau)\,v(z,\,\tau) -  \, |\Omega_j|\, \Phi(x,\,t;\,z_j,\,\tau)\,v(z_j,\,\tau) \right]\,dzd\tau.
\end{equation*}
Following the similar steps as for estimating $\mathcal B_l$, and as the integrands are smooth here, it can be proved that $\mathcal C = o(\varepsilon^{\frac{1}{3}})$ as $\varepsilon \to 0$. Then we have
\begin{equation*}
W(x,\,t) = \sum_{j=1}^{[a^{-1}]} \overline C \, |\Omega_j| \int_0^t \,\Phi(x,\,t;\,z_j,\,\tau)\,\alpha_j(\tau)\,d\tau + \mathcal D + o(\varepsilon^{\frac{1}{3}})
\end{equation*}
with
\begin{equation*}
\mathcal D := -\sum_{j=1}^{[a^{-1}]} \overline C \, |\Omega_j| \int_0^t \,\Phi(x,\,t;\,z_j,\,\tau)\,[\alpha_j(\tau)-v(z_j,\,\tau)]\,d\tau.
\end{equation*}
The term $\mathcal D$ can be estimated as
\begin{eqnarray*}
\mathcal D &=& O\Big(\sum_{j=1}^{[a^{-1}]} |\Omega_j| \int_0^t \vert\Phi(x,\,t;\,z_j,\,\tau)\vert \,\vert \alpha_j(\tau)-v(z_j,\,\tau)\vert\,d\tau\Big)\\
&=& O\Big({a}\sum_{j=1}^{[a^{-1}]} \Big(\int_0^t \vert\Phi(x,\,t;\,z_j,\,\tau)\vert^2\,d\tau\Big)^{1/2}\, \Big(\int^t_0 \vert \alpha_j(\tau)-v(z_j,\,\tau)\vert^2\,d\tau\Big)^{1/2} \Big),
\end{eqnarray*}
and then
\begin{equation*}
\mathcal D=O\left(a\, M^{1/2}a^{-1/6}\right)=O\left(a^{\frac{1}{3}}\right).
\end{equation*}
Hence, we conclude that
\begin{equation}\label{med_final}
W(x,\,t) = u(x,\,t) + O\left(a^{\frac{1}{3}}\right) \quad \mathrm{as}\; a \to 0.
\end{equation}

\bigskip
{\bf Acknowledgement:} The first author was partially supported by the Austrian Science Fund(FWF): P28971-N32. The second author is supported by National Natural Science Foundation of China (No. 11671082) and Qing Lan Project of Jiangsu Province. The work is completed after exchanged visits between RICAM and the School of Mathematics at Southeast University. The authors thank both the two institutions for the arrangements and the friendly atmosphere.

\end{document}